\documentclass[twocolumn]{svjour3}
\usepackage{graphicx}
\usepackage{color}
\usepackage{amsmath,amssymb,calc}
\usepackage{amsfonts}
\newcommand\norm[1]{\left\lVert#1\right\rVert}
\usepackage{mathtools} 
\usepackage{hyperref}
\usepackage[square, numbers, sort&compress]{natbib}   
\usepackage{soul}
\usepackage{subfig}
\usepackage[linesnumbered,ruled]{algorithm2e}

\usepackage[T1]{fontenc}
\let\orgautoref\autoref

\renewcommand{\autoref}[1]{
\def\equationautorefname{Eq.}%
\def\figureautorefname{Fig.}%
\def\sectionautorefname{Sec.}%
\def\subsectionautorefname{Subsec.}%
\def\algorithmautorefname{Alg.}%
\def\subfigureautorefname{Fig.}%
\orgautoref{#1}
}
\makeatletter

\newcommand\Autoref[1]{\@first@ref#1,@}
\def\@throw@dot#1.#2@{#1}
\def\@set@refname#1{
    \edef\@tmp{\getrefbykeydefault{#1}{anchor}{}}%
    \xdef\@tmp{\expandafter\@throw@dot\@tmp.@}%
    \ltx@IfUndefined{\@tmp autorefnameplural}%
         {\def\@refname{\@nameuse{\@tmp autorefname}s}}%
         {\def\@refname{\@nameuse{\@tmp autorefnameplural}}}%
}
\def\@first@ref#1,#2{%
  \ifx#2@\autoref{#1}\let\@nextref\@gobble
  \else%
    \@set@refname{#1}
    \@refname~\ref{#1}
    \let\@nextref\@next@ref
  \fi%
  \@nextref#2%
}
\def\@next@ref#1,#2{%
   \ifx#2@ to~\ref{#1}\let\@nextref\@gobble
   \else, \ref{#1}
   \fi%
   \@nextref#2%
}

\def\tagform@#1{\maketag@@@{\ignorespaces#1\unskip\@@italiccorr}}
\let\orgtheequation\theequation
\def\theequation{(\orgtheequation)}
\makeatother
\journalname{Computational Mechanics}
\begin{document}\sloppy
\title{Graded-material Design based on Phase-field and Topology Optimization}
\author{Massimo Carraturo$^{1,2}$ \and Elisabetta Rocca$^{3,4}$ \and Elena Bonetti$^{4,5}$ \and Dietmar H{\"o}mberg$^{6,7}$ \and Alessandro Reali$^1$ \and Ferdinando Auricchio$^1$
}                     
\authorrunning{Carraturo, M. \textit{et al.}}
%
%
\institute{
{Massimo} {Carraturo} 
\\
massimo.carraturo01@universitadipavia.it
\linebreak
\\
$^1$ Dipartimento di Ingegneria Civile ed Architettura (DICAr), Universit\'a degli Studi Pavia, via Ferrata 3, 27100 Pavia, Italy \\
$^2$ Chair for Computation in Engineering, 
Technical University of Munich, Arcisstr. 21, 80333 Munich, Germany\\
$^3$ Dipartimento di Matematica, Universit\'a degli Studi Pavia, via Ferrata 5, 27100 Pavia, Italy\\
$^4$ IMATI-CNR, via Ferrata 1, 27100 Pavia, Italy\\
$^5$ Dipartimento di Matematica "F.Enriques", Universit\'a degli Studi di Milano,
via Saldini 50, 20133 Milano, Italy\\
$^6$ Weierstrass Institute for Applied Analysis and Stochastics, Mohrenstr. 39, 10117 Berlin, Germany\\
$^7$ Department of Mathematical Sciences, NTNU, Alfred Getz vei 1, 7491 Trondheim, Norway
}
\date{Received: 19 November 2019 / Revised version: \today}
%
\maketitle
\begin{abstract}
In the present work we introduce a novel graded-material design based on phase-field and topology optimization. 
The main novelty of this work comes from the introduction of an additional phase-field variable in the classical single-material phase-field topology optimization algorithm.
This new variable is used to grade the material properties in a continuous fashion.
Two different numerical examples are discussed, in both of them, we perform sensitivity studies to asses the effects of different model parameters onto the resulting structure.
From the presented results we can observe that the proposed algorithm adds additional freedom in the design, exploiting the higher flexibility coming from {additive manufacturing technology}.
\end{abstract}
%
\section{Introduction} \label{sec:introduction}
Structural topology optimization (TO)~\citep{bendsoe_88} is a numerical method which aims, by means of a density function, at optimally distributing a limited amount of material within a volume, representing the initial geometry of a body undergoing specific loads and displacement boundary conditions.
\par 
Structural TO was originally introduced as a discrete formulation where areas of dense material and voids are alternated without any transition region~\citep{bendsoe_83}.
This first approach, also known as the \textit{0-1 topology optimization} problem, leads to many difficulties from both an analytical and a numerical point of view~\citep{sigmund_98}.
\par
A possible alternative approach is based on homogenization methods, where the macroscopic material properties are obtained from microscopic porous material characteristics~\citep{allaire_2004,suzuki_91}.
Such a strategy leads to optimized structures with large \textit{grey-scale} regions of perforated, porous material, which are in general undesired due to their elevated manufacturing complexity, in particular when classical manufacturing processes, such as milling or molding are adopted.
To obtain a clear \textit{black and white} design, \textit{Solid Isotropic Material Penalization} (SIMP) method has been introduced in~\citep{zhou_91}.
The SIMP method consists of penalizing the density region, different from the void or bulk material, by choosing a suitable interpolation scheme for material properties at the macroscopic scale~\citep{bendsoe_99,bendsoe_2003}.
This approach has been successfully employed in many engineering applications beside of structural problems, e.g. fluid analysis~\citep{Gersborg_2005}, fluid- and acoustic-structure interaction~\citep{Yoon_2010,Yoon_2007}, heat conduction~\citep{Gersborg_2006}, multi-physics~\citep{Andreasen_2013}, and composite structures~\citep{Sigmund_1996}.
\par
An alternative to the SIMP method is a TO based on the phase-field method, for the first time introduced by~\citet{bourdin_2003}.
Successively, this method has been employed by~\citet{burger_2006} for stress constrained problems and 
by\\
~\citet{takezawa_2010} in the shape and topology optimization context for minimum compliance and eigenfrequency maximization problems.
More recently,~\citet{penzler_2012} have solved nonlinear elastic problems by means of the phase-field approach,
while~\citet{dede_isogeometric_2012} have been the first to apply this method in the context of isogeometric analysis.
%
Similarly to the SIMP method, phase-field based TO penalizes an approximation of the interface perimeter, such that, by choosing a very small positive penalty term, one can obtain a sharp interface region separating solid materials and voids~\citep{blank_multi-material_2014}.
\par
Inspired by the aforementioned works on phase-field and TO, we aim here at developing an optimization procedure to obtain structures which exploit the possibility of additive manufacturing (AM) technology to distribute material through a body locally varying the material density, leading to the so called functionally graded materials (FGM). 
Numerical models to simulate manufacturing processes of FGM products have been proposed and validated in~\cite{yan_2016,gan_2017a,gan_2017b,yang_2018,wolff_2019} for different AM technologies and materials. All these contributions show how the manufacturing process plays a crucial role in the resulting mechanical properties of FGM structures. Since including the influence of the process parameters within topology optimization schemes is an extremely challenging task (cf.~\citep{allaire_2018} for a first attempt in this direction) we prefer to neglect the influence of AM process on the optimized structure. The presented results have to be thus considered only as an initial starting point towards the final design of optimized FGM structures.

Topology optimization routines to obtain FGM design have been recently investigated by many researchers (see, e.g.~\cite{brackett_2011,cheng_2015,panesar_2018}). Nevertheless, at the best knowledge of the authors a phase-field approach has not been employed to obtain optimized graded-material structures yet. The phase-field method proposed in this contribution allows deriving a rigorous mathematical analysis of the problem for FGM structures and to obtain a clear separation between areas of material and voids.
Moreover, this work aims at investigating by means of 2D examples the sensitivity of the proposed phase-field method to various numerical and physical parameters. On one hand, the choice of considering only plane stress, two-dimensional structures may limit the validity of the presented results while, on the other hand, it allows us to gain a clear insight into the sensitivity of the method w.r.t. the considered parameters with a negligible computational effort.

Even if our approach can potentially include a multi-material case, in this work we consider only a graded-material design, i.e., where a single material is gradually distributed through the body.
The result of such an optimization routine is a structure with graded stiffness values, i.e., a material with stiffness continuously varying within the body, alternating regions of soft material with other regions of stiffer material.
The approach proposed in the present paper reintroduces the typical \textit{grey-scale} regions of early topology optimization methods but within a controlled and numerically stable formulation.
This choice is justified by the fact that modern AM technologies allow grading the density of a body in an almost continuous fashion, varying the amount of distributed material point-by-point during the printing process.
\par
The outline of the work at hand is organized as follows. In~\autoref{sec:SingleMaterial} we recall the formulation for a single-material phase-field based TO.
\autoref{sec:ProblemMultiMat} introduces the novel phase-field approach for graded material structures.
Then,~\autoref{sec:NumericalExamples} discusses two-dimensional numerical examples, carrying out sensitivity studies for different choices of numerical and physical parameters and presenting also a possible solution to obtain a manufactured product from the numerical results.
Finally, in~\autoref{sec:conclusions}, we draw our conclusions on the present work.

\section{Single-material phase-field topology optimization} \label{sec:SingleMaterial}
%
In this section we recall the classical formulation for a phase-field TO of a single-material homogeneous structure, closely following~\citep{BGFS14}.
\subsection{State equations}
\label{ssec:StateEq}
We consider a domain $\Omega\subset\mathbb{R}^{d}$ where material is distributed by means of a scalar phase-field variable $\phi$, representing a material density fraction, hence 
$\phi\in\left[0,1\right]$ with $\phi\equiv 0$ corresponding to voids (i.e., no material) and $\phi\equiv 1$ to bulk material.
Adopting a linear elastic model, the state equations are as follows:
\begin{align}\label{StateEq}
-\text{div}\left(\boldsymbol{\sigma}\right)&=\mathbf{0}\quad\text{in}\quad\Omega \\
\mathbf{u}&=\mathbf{0}\quad\text{on}\quad\Gamma_D\\
\boldsymbol{\sigma}\cdot\mathbf{n}&=\mathbf{g}\quad\text{on}\quad\Gamma_N
\label{StateEqNeumann}
\end{align}
with $\boldsymbol{\sigma}=\boldsymbol{\sigma}(\phi)=\mathbb{C}(\phi)\colon\boldsymbol{\varepsilon}(\mathbf{u})$, $\mathbb{C}=\mathbb{C}(\phi)$ the fourth-order linear material tensor, 
$\mathbf{u}$ the displacement field vector,
$\boldsymbol{\varepsilon}(\mathbf{u})$ the symmetric strain defined as $\boldsymbol{\varepsilon}=\nabla^{S}\mathbf{u}=\left(\nabla\mathbf{u}+\nabla\mathbf{u}^T\right)/2$, 
$\mathbf{g}$ the external load on the boundary $\Gamma_N \subset \partial \Omega$, $\mathbf{n}$ the unit normal vector,
and $\Gamma_D\subset\partial\Omega$, $\mid\Gamma_D\mid\neq0$, the portion of the boundary where homogeneous Dirichlet boundary conditions are applied.
\par
Assuming the material tensor $\mathbb{C}$ to depend on $\phi$, 
the solution of problem~\Autoref{StateEq,StateEqNeumann} depends on the distribution of the scalar field $\phi$ (i.e., $\mathbf{u}=\mathbf{u}(\phi)$). 
We treat the void as a {\sl very soft} material, adopting the following expression for $\mathbb{C}$: 
\begin{equation*}
\mathbb{C}(\phi)=\mathbb{C}_{bulk} \phi^p+  \mathbb{C}_{void} (1-\phi)^p
\end{equation*}
where $\mathbb{C}_{bulk}$ is the positive definite material tensor of the bulk, dense material, $\mathbb{C}_{void}$ is the positive definite material tensor of an idealized very soft material (representing the voids), and $p$ can be any positive value;
for simplicity, we  assume $\mathbb{C}_{void} = \gamma^2 \mathbb{C}_{bulk}$, with $\gamma \ll 1$, while, following~\citep{bendsoe_99}, we set $p=3$.
\par
The weak form of the linear elastic problem~\Autoref{StateEq,StateEqNeumann} can be written as:
\begin{equation}
\int_{\Omega}\boldsymbol{\sigma}(\phi)\colon\boldsymbol{\varepsilon}(\mathbf{v})\text{d}\Omega=\int_{\Gamma_N}\mathbf{g}\cdot \mathbf{v}\text{d}\Gamma.
\label{eq:WeakFormStateEq2}
\end{equation}
with $\mathbf{v}\in\mathcal{H}^1_D(\Omega)$ a virtual displacement field.
Referring to~\citep{BGFS14} we can prove that for any given $\mathbf{g}\in L^2(\Gamma_N)$ and $\phi\in L^\infty(\Omega)$, there exists a unique $\mathbf{u}\in H^1_D(\Omega)$ fulfilling \autoref{eq:WeakFormStateEq2}, with $H^1_D(\Omega):=\{\mathbf{v}\in H^1(\Omega)\,:\, \mathbf{v}=\mathbf{0}\quad \hbox{on }\Gamma_D\}$.
\subsection{Single-material topology optimization as a minimization problem}
\label{ssec:GradientEq}
The goal of our TO process is to properly minimize the compliance of a given structure, by optimally distributing  a limited amount of material. 
\par
To properly minimize the compliance, we introduce an objective functional $\mathcal{J}(\phi,\mathbf{u}(\phi))$ defined as:
\begin{align}
\begin{split}
&\mathcal{J}(\phi,\mathbf{u}(\phi))=
\\
&\int_{\Gamma_N}\mathbf{g}\cdot \mathbf{u}(\phi)\text{d}\Gamma+
\kappa \int_{\Omega}\left[\dfrac{\gamma}{2}\parallel\nabla\phi\parallel^2+\dfrac{1}{\gamma}\psi_0(\phi)\right]\text{d}\Omega
\end{split}
\label{eq:PotentialSingle}
\end{align}
where the first integral represents a measure of the global system compliance, defined as the inverse of the stiffness,
while, assuming  $\kappa > 0$ and a double-well potential function $\psi_0(\phi)=(\phi-\phi^2)^2$ , the second integral is an approximation of the perimeter of the interfaces between regions with $\phi=0$ and $\phi=1$.
In \autoref{eq:PotentialSingle}
$\gamma$ corresponds to the thickness of the diffuse interface, i.e., the region where $0<\phi<1$, the term $\gamma/2\mid\nabla\phi\mid^2$ penalizes jumps between $\phi=0$ and $\phi=1$, while $\psi_0(\phi)/\gamma$ represents the double-well potential function penalizing phases with $\phi$ different from 0 and 1.
We remark that following~\citep{BGFS14} we choose the same scaling parameter $\gamma$ to penalize the sharp interface region and to define the void soft material; this choice is justified by the assumption that when one of the two values goes to zero also the other one has to vanish.
\par
The minimization of the functional in~\autoref{eq:PotentialSingle}
is imposed under the assumption of distributing a limited constant quantity of material inside the domain, 
hence, we introduce the constraint:
\begin{equation*}
\int_{\Omega}\phi\text{d}\Omega = m\mid\Omega \mid
\end{equation*}
with $0 < m \leq 1$ representing a target domain volume fraction.
Clearly, the displacement field $\mathbf{u}(\phi)$ solving the TO problem should also be the solution of the linear elastic problem of~\autoref{eq:WeakFormStateEq2}. 
\par
In conclusion, the minimization problem we aim to solve is the following. 
\\
Problem $(\mathcal{P})$:
\begin{equation*}
\min_{\phi}\quad \mathcal{J}(\phi,\mathbf{u}(\phi))
\label{eq::MinProblemI}
\end{equation*}
such that the following constraints are satisfied:
\begin{align}
&\int_{\Omega}\boldsymbol{\sigma}(\phi)\colon\boldsymbol{\varepsilon}(\mathbf{v})\text{d}\Omega=\int_{\Gamma_N}\mathbf{g}\cdot \mathbf{v}\text{d}\Gamma.
\label{eq:WeakFormStateEq}
\\
& \mathcal{M}(\phi)=\int_{\Omega}\phi\text{d}\Omega - m\mid\Omega\mid=0,
\label{eq:VolumeConstrain}
\end{align}
with $\phi\in H^1(\Omega)$ satisfying the constraint:
\begin{equation}
0\leq\phi\leq 1 \qquad\hbox{a.e. in }\Omega.
\label{eq:PhiConstrain}
\end{equation}
\par
Following the argument by~\citep{BGFS14}, we can prove that the 
minimum constrained problem $(\mathcal{P})$ has at least one solution (cf.~\cite[Thm.~4.1]{BGFS14}).
In particular, to solve problem $(\mathcal{P})$ we introduce the  Lagrangian functional $\mathcal{L}$, defined as:
\begin{equation}
\mathcal{L}(\phi,\mathbf{u},\lambda,\mathbf{p}) = \mathcal{J}(\phi,\mathbf{u}) + \lambda \mathcal{M}(\phi) + \mathcal{S}(\phi,\mathbf{u},\mathbf{p}),
\label{eq:lagrangian}
\end{equation}
where $\lambda$ is the Lagrange multiplier introduced to impose the volume constrain of~\autoref{eq:VolumeConstrain}
and the operator $\mathcal{S}$ is defined as:
\begin{equation*}
\mathcal{S}(\phi,\mathbf{u},\mathbf{p})=\int_{\Omega}\boldsymbol{\sigma}(\phi)\colon\boldsymbol{\varepsilon}(\mathbf{p})\text{d}\Omega -  \int_{\Gamma_N} {\bf g}\cdot {\bf p} \text{d}\Gamma,
\end{equation*}
which we introduce together with the adjoint variable $\mathbf{p}$.
The solution of problem~$(\mathcal{P})$ 
is equivalent to the minimization of~\autoref{eq:lagrangian} 
subjected to constraint in~\mbox{\autoref{eq:PhiConstrain};}
this last problem can be seen as an optimal control problem, with solutions $(\bar{\phi},\bar{\mathbf{u}},\bar{\lambda},\bar{\mathbf{p}})$ that have to satisfy the first order optimality conditions defined by:
\begin{align*}
D_{\mathbf{u}}\mathcal{L}
\left(\bar{\phi},\bar{\mathbf{u}},\bar{\lambda},\bar{\mathbf{p}}\right)
&= 0, \\
D_{\mathbf{p}}\mathcal{L}
\left(\bar{\phi},\bar{\mathbf{u}},\bar{\lambda},\bar{\mathbf{p}}\right)
&= 0, \\
D_{\lambda}\mathcal{L}
\left(\bar{\phi},\bar{\mathbf{u}},\bar{\lambda},\bar{\mathbf{p}}\right)
&= 0, \\
D_{\phi}\mathcal{L}
\left(\bar{\phi},\bar{\mathbf{u}},\bar{\lambda},\bar{\mathbf{p}}\right)\left(\phi-\bar{\phi}\right)
&\geq 0&\forall\phi\in\Phi_{ad},\label{firstOrderOptimalityCondition}
\end{align*}
where $\Phi_{ad}$ is the set of admissible controls defined as follows:
\begin{equation*}
\Phi_{ad}:=\{\phi\in H^1_D(\Omega)\,:\, 0\leq \phi\leq 1 \quad \hbox{a.e. in }\Omega
\}.
\end{equation*}
We also note that for the problem under investigation $D_{\mathbf{p}}\mathcal{L}=D_{\mathbf{u}}\mathcal{L}$, hence the so-called adjoint equation (holding true for every $\mathbf{v}\in H^1_D(\Omega)$):
\begin{equation*}
\int_{\Omega}\mathbb{C}(\phi)\boldsymbol{\varepsilon}(\bar{\mathbf{p}})\colon\boldsymbol{\varepsilon}(\mathbf{v})\text{d}\Omega=\int_{\Gamma_N}\mathbf{g}\cdot \mathbf{v}\text{d}\Gamma,
\label{eq:Adj}
\end{equation*}
is identical to the weak form of the linear elastic problem~\autoref{eq:WeakFormStateEq2}, which implies that $\bar{\mathbf{p}}=\bar{\mathbf{u}}$ .
We refer to~\citep{BGFS14} for the complex analysis of optimality conditions.
\par
To obtain a more compact formulation, we define here the energy density of the system and its derivative w.r.t. the scalar field $\phi$ as:
\begin{equation*}\label{eq::EnergyDensity}
\mathcal{E}(\phi,\mathbf{u})=\boldsymbol{\sigma}(\phi)\colon\boldsymbol{\varepsilon}(\mathbf{u}),
\end{equation*}
and
\begin{equation*}\label{eq::EnergyDensityDerivative}
\dfrac{\partial\mathcal{E}(\phi,\mathbf{u})}{\partial\phi}=\dfrac{\partial\boldsymbol{\sigma}(\phi)}{\partial\phi}\colon\boldsymbol{\varepsilon}(\mathbf{u}),
\end{equation*}
where 
\begin{equation*}
\dfrac{\partial\boldsymbol{\sigma}(\phi)}{\partial\phi} = \dfrac{\partial\mathbb{C}(\phi)}{\partial\phi}\colon\boldsymbol{\varepsilon}(\mathbf{u}).
\end{equation*}
\par
To discretize our continuous problem we employ a gradient flow dynamics, namely Allen-Cahn gradient flow~\citep{allen_microscopic_1979}, a steepest descent pseudo-time stepping method with a time-step increment $\tau$.
Thus the optimal control problem $(\mathcal{P})$ can be now rewritten as follows:
\begin{align}\label{eq::FunctionalForm}
&D_{\mathbf{u}}\mathcal{L}\:\mathbf{v} =0,\\
&D_{\lambda}\mathcal{L}\:v_{\lambda}=\mathcal{M}v_{\lambda}=0,\\
&\dfrac{\gamma}{\tau}\int_{\Omega}(\phi_{n+1}-\phi_n)v_{\phi}\text{d}\mathbf{x}=-D_\phi \mathcal{L}\:v_{\phi},\label{eq::AllenCahnFunc}
\end{align}
where
\begin{equation*}
D_{\phi}\mathcal{L} = \dfrac{\partial\mathcal{J}}{\partial\phi} + \lambda\dfrac{\partial\mathcal{M}}{\partial\phi} + \dfrac{\partial\mathcal{S}}{\partial\phi},
\label{eq::LagrangianDerPhi}
\end{equation*} 
with $v_{\lambda}\in\mathbb{R}$ and $v_{\phi} \in \Phi_{ad}$.
\par
The problem defined in~\Autoref{eq::FunctionalForm, eq::AllenCahnFunc} can be written in the following weak extended formulation:
\begin{align}
&\int_{\Omega}\boldsymbol{\sigma}(\phi)\colon\boldsymbol{\varepsilon}(\mathbf{v})\text{d}\Omega=\int_{\Gamma_N}\mathbf{g}\cdot \mathbf{v}\text{d}\Gamma,
\label{eq:WeakFormStateEq3}
\\
&\int_{\Omega}v_{\lambda}(\phi_{n+1}-m)\text{d}\Omega
 = 0,
\label{eq:GradientEqLambda}
\\
&\dfrac{\gamma}{\tau}\int_{\Omega}(\phi_{n+1}-\phi_n)v_{\phi}\text{d}\Omega +
\kappa\gamma\int_{\Omega}\nabla\phi_{n+1}\cdot\nabla v_{\phi}\text{d}\Omega 
\notag\\
&\phantom{{}=2}+\lambda \int_{\Omega}v_{\phi}\text{d}\Omega - \int_{\Omega}v_{\phi}\dfrac{\partial\mathcal{E}(\phi_n,\mathbf{u}_n)}{\partial\phi}\text{d}\Omega
\notag\\
&\phantom{{}=2}+\dfrac{\kappa}{\gamma}\int_{\Omega}\dfrac{\partial\psi_0(\phi_{n})}{\partial\phi}v_{\phi}\text{d}\Omega=0.
\label{eq:GradientEqPhi}
\end{align}
\subsection{Single-material finite element formulation}
\label{ssec:FEM}
We derive here a finite element approximation of the phase-field TO problem defined in~\Autoref{eq:WeakFormStateEq3,eq:GradientEqPhi}.
To this end we discretize the physical domain $\Omega$ using two different meshes $\mathcal{Q}_u$, $\mathcal{Q}_{\phi}$ corresponding to the field variables $\mathbf{u}$ and $\phi$ and their variations $\mathbf{v}$ and $v_{\phi}$. The Lagrange multiplier $\lambda$ used to constrain the volume is applied using a constant scalar value on $\Omega$.
On each mesh, we interpolate the nodal values of the field variables and their variations by means of piecewise linear basis functions, such that:
\begin{align*}
& \mathbf{u}\approx\mathbf{N_u}\tilde{\mathbf{u}}, & \mathbf{v}\approx\mathbf{N_u}\tilde{\mathbf{v}}, \\
& \phi\approx\mathbf{N}_{\phi}\tilde{\boldsymbol{\phi}}, & v_{\phi}\approx\mathbf{N}_{\phi}\tilde{\mathbf{v}}_{\phi}, \\
\end{align*}
\par
Introducing the proposed discretization in~\Autoref{eq:WeakFormStateEq3,eq:GradientEqPhi} the discrete version of the optimal control problem becomes:
\begin{equation}
\begin{split}
\renewcommand\arraystretch{1.3}
\dfrac{1}{\tau}
\left[
\begin{matrix}
 \mathbf{0} & \mathbf{0} & \mathbf{M}^{\phi\lambda} \\
 \mathbf{0} &  \mathbf{M}^{\phi\phi} & \mathbf{0} \\
 \mathbf{0} &  \mathbf{M}^{\lambda\phi} & \mathbf{0}
\end{matrix}
\right]
\left[
\begin{array}{c}
\tilde{\mathbf{u}}  \\
\tilde{\boldsymbol{\phi}} \\
\tilde{\lambda}
\end{array}
\right]
+
\left[
\begin{matrix}
\mathbf{K}^{\mathbf{u}\mathbf{u}} & \mathbf{0}  &  \mathbf{0} \\
\mathbf{0} & \mathbf{K}^{\phi\phi} & \mathbf{0} \\
\mathbf{0} & \mathbf{0}  &  \mathbf{0} \\
\end{matrix}
\right]
\left[
\begin{array}{c}
\tilde{\mathbf{u}}  \\
\tilde{\boldsymbol{\phi}}\\
\tilde{\lambda}
\end{array}
\right]
\\
=
\left[
\begin{array}{c}
\mathbf{f}  \\
\mathbf{q}^{\phi}+\mathbf{q}^{s}+\mathbf{q}^{\psi} \\
q^{\lambda}
\end{array}
\right]
\end{split}
\label{eq:DiscreteSystem}
\end{equation}
with the matrix and vector terms defined as follows:
\begin{align*}
&\mathbf{K}^{\mathbf{u}\mathbf{u}} =  \int_{\Omega}\nabla\mathbf{N_u}^T\mathbb{C}\nabla\mathbf{N_u}\text{d}\Omega,
\\
& \mathbf{M}^{\phi\phi}=\gamma\int_{\Omega}\mathbf{N}^T_{\phi}\mathbf{N}_{\phi}\text{d}\Omega,
\\
& \mathbf{K}^{\phi\phi}=\kappa\gamma\int_{\Omega}\nabla\mathbf{N}_{\phi}^T\nabla\mathbf{N}_{\phi}\text{d}\Omega,
\\
& \mathbf{M}^{\lambda\phi}=\tau\int_{\Omega}\mathbf{N}_{\phi}\text{d}\Omega=\left(\mathbf{M}^{\phi\lambda}\right)^T,
\\
&\mathbf{f} =  \int_{\Gamma_N}\mathbf{N_u}^T\mathbf{g}\text{d}\Gamma,
\\
& \mathbf{q}^{\phi}=\dfrac{\gamma}{\tau}\int_{\Omega}\left(\mathbf{N}^T_{\phi}\mathbf{N}_{\phi}\right)\boldsymbol{\tilde{\phi}}_n\text{d}\Omega = \mathbf{M}^{\phi\phi}\boldsymbol{\tilde{\phi}}_n,
\\
& q^{\lambda}=\int_{\Omega}m \text{d}\Omega,
\\
& \mathbf{q}^{s}=\int_{\Omega}\mathbf{N}_{\phi}^T \dfrac{\partial\mathcal{E}(\boldsymbol{\tilde{\phi}}_n,\mathbf{\tilde{u}}_n)}{\partial\phi} \text{d}\Omega,
\\
& \mathbf{q}^{\psi}=  -\dfrac{\kappa}{\gamma} \int_{\Omega}\mathbf{N}_{\phi}^T \dfrac{\partial\psi_0(\boldsymbol{\tilde{\phi}}_n)}{\partial\phi} \text{d}\Omega.
\end{align*}
The discrete linear system in~\autoref{eq:DiscreteSystem} can be solved using a staggered approach, i.e. solving first the state equation system:
\begin{equation}
  \mathbf{K^{uu}}\mathbf{\tilde{u}} =
  \mathbf{f},
  \label{eq:DiscreteStateSystem}
\end{equation}
and then the discretized optimization problem:
\begin{equation}
\renewcommand\arraystretch{1.3}
\dfrac{1}{\tau}
\left[
\begin{array}{c c}
\begin{matrix}
\mathbf{M}^{\phi\phi}+\tau\mathbf{K}^{\phi\phi} & \mathbf{M}^{\phi\lambda} \\
\mathbf{M}^{\lambda\phi} & \mathbf{0} 
\end{matrix}
\end{array}
\right]
\left[
\begin{array}{c}
\tilde{\boldsymbol{\phi}}\\
\tilde{\lambda} 
\end{array}
\right]
=
\left[
\begin{array}{c}
\mathbf{q}^{\phi}+\mathbf{q}^{s}+\mathbf{q}^{\psi} \\
q^{\lambda}
\end{array}
\right].
  \label{eq:DiscreteGradientEquationSystem}
\end{equation}
\par
Following this approach~\autoref{optAlgorithm} solves the problem in two separate steps:
first it solves the state equation system~\autoref{eq:DiscreteStateSystem} to get the solution vector $\mathbf{\tilde{u}}_{n+1}$ (line 3),
secondly, the linear system of~\autoref{eq:DiscreteGradientEquationSystem} is solved to obtain the phase-field vector $\boldsymbol{\tilde{\phi}^\ast}_{n+1}$ and the Lagrange multiplier vector $\tilde{\lambda}_{n+1}$~\mbox{(line 4).}
Finally, the vector $\boldsymbol{\tilde{\phi}^\ast}_{n+1}$  is projected within the interval $\left[0,1\right]$ to obtain the phase-field solution vector $\boldsymbol{\tilde{\phi}}_{n+1}$ fulfilling the constraints (line 5).
We use the increment $\Delta_\phi$ based on the $L^2$-norm and defined as:
\begin{equation}
\Delta_\phi=\dfrac{\norm{\boldsymbol{\tilde{\phi}}_{n+1}-\boldsymbol{\tilde{\phi}}_n}_{L^2}}{\norm{\boldsymbol{\tilde{\phi}}_n}_{L^2}},
\end{equation}
as a criterion to assert the convergence of the algorithm, which is otherwise stopped when user defined maximum number of iterations $max_{iter}$ is reached.
\begin{algorithm}
\DontPrintSemicolon
\SetKwInOut{Input}{input}\SetKwInOut{Output}{output}
\Input{$\mathcal{Q}_u$, $\mathcal{Q}_{\phi}$,$\boldsymbol{\tilde{\phi}}_0$.}
\Output{Optimal topology}
$\boldsymbol{\tilde{\phi}}_n\gets \boldsymbol{\tilde{\phi}}_0$\;
\While {$\Delta_\phi\geq tol$ and $n\leq max_{iter}$}
{
		{$\mathbf{\tilde{u}}_{n+1}\gets$\texttt{solve}\eqref{eq:DiscreteStateSystem}\;}
		{$(\tilde{\lambda}_{n+1},\boldsymbol{\tilde{\phi}}^\ast_{n+1})\gets$\texttt{solve}\eqref{eq:DiscreteGradientEquationSystem}\;}
		{$\boldsymbol{\tilde{\phi}}_{n+1}\gets$\texttt{rescale}$\left(\boldsymbol{\tilde{\phi}}^\ast_{n+1}\right)$ to $\left[0,1\right]$\;}
		{\texttt{update}($\Delta_\phi$)\;}
		{$\boldsymbol{\tilde{\phi}}_n\gets \boldsymbol{\tilde{\phi}}_{n+1}$\;}
}
\caption{Single-material optimization algorithm}\label{optAlgorithm}
\end{algorithm}
\section{Graded-material phase-field topology optimization} \label{sec:ProblemMultiMat}
In the following section, we extend the previously presented formulation of TO to the case of a graded material definition. We refer to this approach as graded-material phase-field topology optimization.
The mathematical analysis of the corresponding optimization problem will be the subject of a forthcoming paper~\citep{Rocca_2018}.
\subsection{State equation}
\label{ssec:MMTensor}
We now consider the case of an inhomogeneous material distribution; in particular, we assume that the material elastic fourth-order tensor $\mathbb{C}$ can vary linearly through a material grading scalar variable $\chi\in\left[0,\phi\right]$, with $\phi\in\left[0,1\right]$, such that:

\begin{equation*}
\mathbb{C}(\chi)=\mathbb{C}_{bulk}\chi+\dfrac{1}{\beta}\mathbb{C}_{bulk}(\phi-\chi),
\end{equation*} 
with $0 < \beta \leq 1$ a so-called softening factor, used to define the soft material tensor as a fraction of the bulk, fully dense material tensor $\mathbb{C}_{bulk}$.
In such a way, the stiffness of the body can continuously vary from a full dense material ($\chi=\phi$) to a softer one ($\chi=0$). Therefore, the definition of the fourth-order material tensor $\mathbb{C}(\phi,\chi)$, previously defined in~\autoref{sec:SingleMaterial}, can now be modified as:
\begin{equation}\label{eq:MMTensor}
\mathbb{C}(\phi,\chi)=\mathbb{C}(\chi)\phi^p+\gamma_{\phi}^2\mathbb{C}(\chi)(1-\phi)^p,
\end{equation}
where $0<\gamma_{\phi}\ll 1$ and again we choose a penalty parameter $p=3$.

We would like to remark here that FGM structures are intrinsically heterogeneous but, as recently demonstrated by~\citet{cheng_2019}, an asymptotic homogenization method can be effectively employed. In this work, we assume a simple linearly interpolated homogenized material tensor but a more complex homogenization could be directly introduced within the proposed numerical scheme.

The definition of the material tensor in~\autoref{eq:MMTensor} leads to an optimized structure where, as in the single-material case, the perimeter of the body is defined by the sharp interface of the phase-field variable $\phi$, while the stiffness of the material continuously varies within the structure, following the distribution of the material grading variable $\chi$.
Hence, the graded-material weak form of the linear elastic problem of~\autoref{eq:WeakFormStateEq2} reads: 
\begin{equation}\label{eq:StateEqMM}
\int_{\Omega}\boldsymbol{\sigma}(\phi,\chi)\colon\boldsymbol{\varepsilon}(\mathbf{v})\text{d}\Omega=\int_{\Gamma_N}\mathbf{g}\cdot \mathbf{v}\text{d}\Gamma
\end{equation}
with the virtual displacement $\mathbf{v}\in H^1_D(\Omega)$ and where $\boldsymbol{\sigma}(\phi,\chi)=\mathbb{C}(\phi,\chi)\colon\boldsymbol{\varepsilon}(\mathbf{u})$.
%
\subsection{Graded-material topology optimization as a minimization problem}
\label{ssec:MMGradEq}
We want now to define an objective functional which optimizes a structure with an inhomogeneous material distribution.
This new graded-material objective functional $\mathcal{J}^M(\phi,\chi,\mathbf{u}(\phi,\chi))$ can be defined as:
\begin{align*}
\begin{split}
&\mathcal{J}^M(\phi,\chi,\mathbf{u}(\phi,\chi))=\int_{\Gamma_N}\mathbf{g}\cdot \mathbf{u}(\phi,\chi)\text{d}\Gamma+
\\
&\kappa_{\phi} \int_{\Omega}\left[\dfrac{\gamma_{\phi}}{2}\mid\nabla\phi\mid^2+
\dfrac{1}{\gamma_{\phi}}\psi_0(\phi)\right]\text{d}\Omega+
\kappa_{\chi}\int_{\Omega}\dfrac{\gamma_{\chi}}{2}\mid\nabla\chi\mid^2\text{d}\Omega,
\end{split}
\end{align*}
with $\kappa_{\phi},\kappa_{\chi}>0$ and $\gamma_{\chi}>0$, and where the first two integrals are the same of the objective functional in~\autoref{eq:PotentialSingle}, while the additional integral term~$\gamma_{\chi}/2\left(\mid\nabla\chi\mid\right)^2$ is introduced to penalize the gradient of the scalar field $\chi$.
\par
Following the same approach described for the single-material case, the global graded-material minimization problem can be now written as follows:
\begin{align*}
& \min_{\phi,\chi}\quad \mathcal{J}^M(\phi,\chi,\mathbf{u}(\phi,\chi)),
\\
\text{such that:}
\\
&\int_{\Omega}\boldsymbol{\sigma}(\phi,\chi)\colon\boldsymbol{\varepsilon}(\mathbf{v})\text{d}\Omega=\int_{\Gamma_N}\mathbf{g}\cdot \mathbf{v}\text{d}\Gamma,
\\
& \mathcal{M}(\phi)=\int_{\Omega}\phi\text{d}\Omega- m\mid\Omega\mid=0,
\end{align*}
where $\phi,\chi\in H^1(\Omega)$, under the constraint
\begin{equation*}
0\leq\phi\leq 1 \qquad\hbox{a.e. in }\Omega,
\label{eq:MMPhiConstrain}
\end{equation*}
and the additional constraint on $\chi$:
\begin{equation*}\label{eq:ChiConstrain}
0 \leq \chi \leq \phi \quad \hbox{a.e. in }\Omega.
\end{equation*}
We can now define the graded-material Lagrangian $\mathcal{L}^M$ as:
\begin{align*}
\mathcal{L}^M
=
\mathcal{J}^M
+ \lambda\mathcal{M} + \mathcal{S}^M,
\end{align*} 
explicitly written as:
\begin{align*}
\mathcal{L}^M(\phi,\chi,\mathbf{u},\lambda,\mathbf{p}) =
\mathcal{J}^M(\phi,\chi,\mathbf{u}) + \lambda\mathcal{M}(\phi) + \mathcal{S}^M(\phi,\chi,\mathbf{u},\mathbf{p}),
\label{eq:lagrangianMM}
\end{align*} 
where, the operator $\mathcal{S}^M$ for the graded-material formulation is calculated as:
\begin{equation*}
\mathcal{S}^M(\phi,\chi,\mathbf{u},\mathbf{p}) = \int_{\Omega}\boldsymbol{\sigma}(\phi,\chi)\colon\boldsymbol{\varepsilon}(\mathbf{p})\text{d}\Omega - \int_{\Gamma_N} {\bf g}\cdot {\bf p} \text{d}\Gamma.
\end{equation*} 
Analogously to the previously introduced set of admissible controls $\Phi_{ad}$ for the phase-field variable $\phi$, we define now the set of admissible controls $\Xi_{ad}$ for the grading variable $\chi$ as:
\begin{equation*}
\Xi_{ad}:=\{\chi\in H^1(\Omega)\,:\, 0\leq \chi\leq \phi \quad \hbox{a.e. in }\Omega\}.
\end{equation*}
\par
Clearly, also in the graded-material case, we want that the optimal control solutions $\bar{\phi}$ and $\bar{\chi}$ have to satisfy the first order necessary optimality conditions, which can be derived as:
\begin{equation*}
D_{\phi}\mathcal{L}^M(\bar{\phi},\bar{\chi},\bar{\mathbf{u}},\bar{\lambda},\bar{\mathbf{p}})\left(\phi-\bar{\phi}\right)\geq 0\quad\forall\phi\in\Phi_{ad}
\label{OptimalityMMI}
\end{equation*}
and
\begin{equation*}
D_{\chi}\mathcal{L}^M\left(\bar{\phi},\bar{\chi},\bar{\mathbf{u}},\bar{\lambda},\bar{\mathbf{p}}\right)\left(\chi-\bar{\chi}\right)\geq 0\quad\forall\chi\in\Xi_{ad},
\label{OptimalityMMII}
\end{equation*}
where
$\bar{\mathbf{u}}$
and
$\bar{\mathbf{p}}$ are solutions of the graded-material state equation~\autoref{eq:StateEqMM} and of the corresponding adjoint problem, respectively.
As in the previous case, the displacement field $\mathbf{u}$ is self-adjoint and hence we have $\bar{\mathbf{p}}=\bar{\mathbf{u}}$.
For a complete analysis of necessary first order optimality conditions we refer to the forthcoming paper~\citep{Rocca_2018}.
\par
Analogously to the single-material case, we can define the energy density of the system and its derivatives w.r.t. both the scalar field $\phi$ and the material grading variable $\chi$ as:
\begin{equation*}\label{eq::EnergyDensityMM}
\mathcal{E}^M(\phi,\chi,\mathbf{u})=\boldsymbol{\sigma}(\phi,\chi)\colon\boldsymbol{\varepsilon}(\mathbf{u}),
\end{equation*}
\begin{align*}
\begin{split}
&\dfrac{\partial\mathcal{E}^M(\phi,\chi,\mathbf{u})}{\partial\phi}=\dfrac{\partial\boldsymbol{\sigma}(\phi,\chi)}{\partial\phi}\colon\boldsymbol{\varepsilon}(\mathbf{u}) =
\\
& \left[3\mathbb{C}(\chi)\phi^2+3\gamma_{\phi}^2\mathbb{C}(\chi)(1-\phi)^2\right]\boldsymbol{\varepsilon}(\mathbf{u})\colon\boldsymbol{\varepsilon}(\mathbf{u})
\end{split}
\end{align*}
and\\
\begin{align*}
\begin{split}
&\dfrac{\partial\mathcal{E}^M(\phi,\chi,\mathbf{u})}{\partial\chi}=\dfrac{\partial\boldsymbol{\sigma}(\phi,\chi)}{\partial\chi}\colon\boldsymbol{\varepsilon}(\mathbf{u}) = \\
&\left[\left(\mathbb{C}_{bulk}-\dfrac{1}{\beta}\mathbb{C}_{bulk}\right)\left(\phi^3+\gamma_{\phi}^2(1-\phi)^3\right)\right]\boldsymbol{\varepsilon}(\mathbf{u})\colon\boldsymbol{\varepsilon}(\mathbf{u}).
\end{split}
\end{align*}
\par
The optimal control problem can be solved as in the single-material case by means of the Allen-Cahn gradient flow, leading to the following set of equations:
\begin{multline}
\dfrac{\gamma_{\phi}}{\tau}\int_{\Omega}(\phi_{n+1}-\phi_{n})v_{\phi}\text{d}\Omega +
\kappa_{\phi}\gamma_{\phi}\int_{\Omega}\nabla\phi\cdot\nabla v_{\phi}\text{d}\Omega +
\\
\int_{\Omega}v_{\phi}\lambda\text{d}\Omega 
- \int_{\Omega}v_{\phi}\dfrac{\partial\mathcal{E}^M(\phi_n,\chi_n,\mathbf{u}_n)}{\partial\phi}\text{d}\Omega+
\\
\dfrac{\kappa_{\phi}}{\gamma_{\phi}}\int_{\Omega}\dfrac{\partial\psi_0(\phi_n)}{\partial\phi}v_{\phi}\text{d}\Omega = 0,
 \label{eq:GradientEqMultiPhi}
\end{multline}
\begin{multline}
\dfrac{\gamma_{\chi}}{\tau}\int_{\Omega}(\chi_{n+1}-\chi_n)v_{\chi}\text{d}\Omega +
\kappa_{\chi}\gamma_{\chi}\int_{\Omega}\nabla\chi\cdot\nabla v_{\chi}\text{d}\Omega -
\\
\int_{\Omega}v_{\chi} \dfrac{\partial\mathcal{E}^M(\phi_n,\chi_n,\mathbf{u}_n)}{\partial\chi} \text{d}\Omega = 0,
 \label{eq:GradientEqMultiChi}
\end{multline}
to be solved under the volume constraint
\begin{equation}
\int_{\Omega}v_{\lambda}(\phi-m)\text{d}\Omega = 0.
 \label{eq:GradientEqMultiLambda}
\end{equation}
In order to estimate the total amount of material in the structure, we define a material fraction index $m_{\chi}$ as:
\begin{equation*}
m_{\chi} = \dfrac{1}{\mid\Omega\mid}\int_{\Omega}\chi\text{d}\Omega,
\end{equation*}
which can be considered as a measure of the global amount of material used to print the structure.
The equivalent material fraction index for the single-material case $m_{\phi}$ is equal to the volume fraction $m$, such that:
\begin{equation*}
m_{\phi} = m = \dfrac{1}{\mid\Omega\mid}\int_{\Omega}\phi\text{d}\Omega.
\end{equation*}
\subsection{Graded-material finite element formulation}
\label{ssec:MultiMaterialFEM}
We aim now at obtaining a discrete formulation for the graded-material phase-field topology optimization problem.
To this end, the displacement field $\mathbf{u}$, the phase-field variable $\phi$, the Lagrange multiplier $\lambda$ and their corresponding variations are approximated using the same discretization already defined in~\mbox{\autoref{ssec:FEM}}.
Additionally, we need to discretize the material grading variable $\chi$ on the domain~\mbox{$\Omega$;} such a discretization is obtained introducing an additional mesh $\mathcal{Q}_{\chi}$, such that the material grading variable $\chi$ and its variation $v_{\chi}$ can be written as:
\begin{align*}
& \mathbf{\chi}\approx\mathbf{N}_{\chi}\tilde{\boldsymbol{\chi}} & \text{and}\qquad\qquad\qquad & v_{\chi}\approx\mathbf{N}_{v_{\chi}}\tilde{\mathbf{v}}_{\chi},
\end{align*}
where
$\mathbf{N}_{\chi}$ and $\mathbf{N}_{v_{\chi}}$ are the piecewise linear shape functions which interpolate the nodal degrees of freedoms $\tilde{\boldsymbol{\chi}}$ and $\tilde{\mathbf{v}}_{\chi}$, respectively.
\par
The discrete form of~\Autoref{eq:GradientEqMultiPhi, eq:GradientEqMultiLambda} can thus be written in a compact notation as:
\begin{multline}
\dfrac{1}{\tau}
\begin{bmatrix}
\mathbf{0}  & \mathbf{0}  & \mathbf{0} & \mathbf{0}  \\
\mathbf{0}  & \mathbf{M}^{\phi\phi} & \mathbf{0}  & \mathbf{M}^{\phi\lambda} \\
\mathbf{0}  & \mathbf{0} & \mathbf{M}^{\chi\chi} & \mathbf{0} \\
\mathbf{0}  & \mathbf{M}^{\lambda\phi} & \mathbf{0} & \mathbf{0}
\end{bmatrix}
\begin{bmatrix}
\tilde{\mathbf{u}}  \\
\tilde{\boldsymbol{\phi}}  \\
\tilde{\boldsymbol{\chi}} \\
\tilde{\lambda}
\end{bmatrix}
+
\begin{bmatrix}
\mathbf{K}^{\mathbf{u}\mathbf{u}} & \mathbf{0}  & \mathbf{0} & \mathbf{0} \\
\mathbf{0} & \mathbf{K}^{\phi\phi} & \mathbf{0}  & \mathbf{0} \\
\mathbf{0} & \mathbf{0} &\mathbf{K}^{\chi\chi} & \mathbf{0} \\
\mathbf{0} & \mathbf{0} & \mathbf{0} & \mathbf{0}
\end{bmatrix}
\begin{bmatrix}
\tilde{\mathbf{u}}  \\
\tilde{\boldsymbol{\phi}}  \\
\tilde{\boldsymbol{\chi}} \\
\tilde{\lambda}
\end{bmatrix}
\\
 =
\begin{bmatrix}
\mathbf{f} \\
\mathbf{q}^{\phi} + \mathbf{q}^{s\prime} + \mathbf{q}^{\psi} \\
\mathbf{q}^{\chi} + \mathbf{q}^{t}\\
q^{\lambda}
\end{bmatrix},
\label{discreteMultiMaterialSystem}
\end{multline}
where the newly defined matrix and vector terms are:
\begin{align*}
& \mathbf{M}^{\chi\chi}=\gamma_{\chi}\int_{\Omega}\mathbf{N}^T_{\chi}\mathbf{N}_{\chi}\text{d}\Omega,\\
& \mathbf{K}^{\chi\chi}=\kappa_{\chi}\gamma_{\chi}\int_{\Omega}\nabla\mathbf{N}_{\chi}^T\nabla\mathbf{N}_{\chi}\text{d}\Omega,\\
& \mathbf{q}^{\chi}=\dfrac{\gamma_{\chi}}{\tau}\int_{\Omega}\mathbf{N}^T_{\chi}\mathbf{N}_{\chi}\boldsymbol{\tilde{\chi}}_n\text{d}\Omega,\\
\begin{split}
& \mathbf{q}^{s\prime}= \int_{\Omega}\mathbf{N}_{\phi}^T \dfrac{\partial\mathcal{E}^M(\boldsymbol{\tilde{\phi}}_n,\boldsymbol{\tilde{\chi}}_n,\mathbf{\tilde{u}}_n)}{\partial\phi}\text{d}\Omega,
\end{split}
\\
& \mathbf{q}^{t}= \int_{\Omega}\mathbf{N}_{\phi}^T \dfrac{\partial\mathcal{E}^M(\boldsymbol{\tilde{\phi}}_n,\boldsymbol{\tilde{\chi}}_n,\mathbf{\tilde{u}}_n)}{\partial\chi}\text{d}\Omega.
\end{align*}
\autoref{optAlgorithmMultiMaterial} describes the iterative procedure to obtain the graded-material optimized structure discussed so far.
The adopted solution scheme is very similar to~\autoref{optAlgorithm} but in this case, we have to solve at each iteration the graded-material linear system defined in~\autoref{discreteMultiMaterialSystem} to obtain the phase-field solution vector $\boldsymbol{\tilde{\phi}}_{n+1}$ and the grading scalar variable vector $\boldsymbol{\tilde{\chi}}_{n+1}$.
As in the single material case, the system can be solved following a staggered scheme, since~\autoref{discreteMultiMaterialSystem} can be split into two separate systems as follows:
\begin{equation}
  \mathbf{K^{uu}}\mathbf{\tilde{u}} =
  \mathbf{f},
  \label{eq:MMDiscreteStateSystem}
\end{equation}
and
\begin{equation}
\renewcommand\arraystretch{1.3}
\dfrac{1}{\tau}
\left[
\begin{array}{c c c}
\begin{matrix}
\mathbf{M}^{\phi\phi}+\tau\mathbf{K}^{\phi\phi} & \mathbf{0}  & \mathbf{M}^{\phi\lambda}\\
\mathbf{0}  & \mathbf{M}^{\chi\chi}+\tau\mathbf{K}^{\chi\chi} & \mathbf{0} \\
\mathbf{M}^{\lambda\phi} & \mathbf{0} & \mathbf{0} \\
\end{matrix}
\end{array}
\right]
\left[
\begin{array}{c}
\tilde{\boldsymbol{\phi}} \\
\tilde{\boldsymbol{\chi}} \\
\tilde{\lambda}\\
\end{array}
\right]
=
\left[
\begin{array}{c}
\mathbf{q}^{\phi}+\mathbf{q}^{s'}+\mathbf{q}^{\psi} \\
\mathbf{q}^{\chi}+\mathbf{q}^{t}\\
q^{\lambda}
\end{array}
\right].
  \label{eq:MMDiscreteGradientEquationSystem}
\end{equation}
Finally, we use the relative increment of $\chi$ in the $L^2$-norm, defined as:
\begin{equation}
\Delta_\chi=\dfrac{\norm{\boldsymbol{\tilde{\chi}}_{n+1}-\boldsymbol{\tilde{\chi}}_n}_{L^2}}{\norm{\boldsymbol{\tilde{\chi}}_n}_{L^2}}.
\end{equation}
as an additional criteria to check the convergence of~\autoref{optAlgorithmMultiMaterial}.
\begin{algorithm}
\DontPrintSemicolon
\SetKwInOut{Input}{input}\SetKwInOut{Output}{output}
\caption{Graded-material optimization algorithm}\label{optAlgorithmMultiMaterial}
\Input{$\mathcal{Q}_u$, $\mathcal{Q}_{\phi}$, $\mathcal{Q}_{\chi}$,$\boldsymbol{\phi}_0$, $\boldsymbol{\chi}_0$}
\Output{Optimal topology}
$\boldsymbol{\phi}_n\gets\boldsymbol{\phi}_0$\;
$\boldsymbol{\chi}_n\gets\boldsymbol{\chi}_0$\;
\While{($\Delta_\phi\geq tol$ or $\Delta_\chi\geq tol$) and $n\leq max_{iter}$}
{
		$\mathbf{\tilde{u}}_{n+1}\gets$ \texttt{solve}\eqref{eq:MMDiscreteStateSystem}\;
		$(\boldsymbol{\tilde{\phi}}^\ast_{n+1},\boldsymbol{\tilde{\chi}}^\ast_{n+1},
		\tilde{\lambda}_{n+1})\gets$ \texttt{solve}\eqref{eq:MMDiscreteGradientEquationSystem}\;
		$\boldsymbol{\tilde{\phi}}_{n+1}\gets$ \texttt{rescale} $\left(\boldsymbol{\tilde{\phi}}^\ast_{n+1}\right)$ to $\left[0,1\right]$\;
				$\boldsymbol{\tilde{\chi}}_{n+1}\gets$ \texttt{rescale} $\left(\boldsymbol{\tilde{\chi}}^\ast_{n+1}\right)$ to $\left[0,\phi\right]$\;
  \texttt{update}($\Delta_\phi$ and $\Delta_\chi$)\;
	    $\boldsymbol{\phi}_n\gets\boldsymbol{\phi}_{n+1}$ \;
		$\boldsymbol{\chi}_n\gets\boldsymbol{\chi}_{n+1}$ \;
}
\end{algorithm}
\section{Numerical Examples} \label{sec:NumericalExamples}
In this section two numerical examples are presented:
in the first one, we consider a cantilever beam structure while in the second one we study a simply-supported beam structure.
For each example, we decided to run two sensitivity studies w.r.t. a numerical and a physical problem parameter to asses the robustness of~\autoref{optAlgorithmMultiMaterial} to these changes.
In the cantilever beam example, we discuss two sensitivity studies. Firstly, we vary the values of the graded-material interface parameter $\gamma_{\chi}$ (i.e., the parameter which represents the thickness of the material grading variable interface) and, secondly, we change the slenderness of the structure.
In the simply-supported beam example, again we perform two studies. In the first one, we use different values for the softening factor $\beta$ while in the second one we increase the load acting on the structure. 
Finally, in~\autoref{ssec:NumToAM} we present a possible solution to obtain an AM product from our numerical results.
\subsection{Cantilever beam}
\label{ssec:CantBeam}
\begin{figure}
\centering
\includegraphics[width=0.45\textwidth]{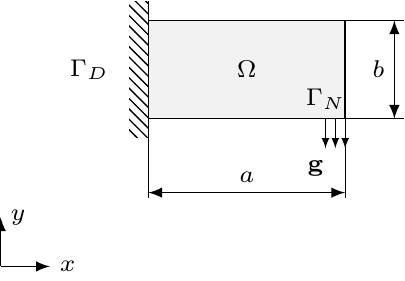}
\caption{Cantilever beam: Initial configuration and problem domain.\label{fig:ProblemSetup}}
\end{figure}

We consider the cantilever beam problem depicted in~\autoref{fig:ProblemSetup}, with dimensions $a=2$mm and $b=1$mm and a traction force $\mathbf{g}=(0,-600)$N/mm applied at the right-end of the lower edge of the structure, while the left edge is fixed.
We assume the initial material being a dense isotropic material, i.e., $\mathbb{C}_{bulk}=(\lambda+2\mu)\mathbf{1}\otimes\mathbf{1}+2\mu\textbf{I}$, where the Lame's parameters $\lambda$ and $\mu$ can be expressed in terms of the Young modulus $E$ and the Poisson coefficient $\nu$ as follows:
\begin{equation}
\lambda=\dfrac{E\nu}{(1+\nu)(1-2\nu)},
\end{equation}
and
\begin{equation}
\mu=\dfrac{E}{2(1+\nu)}.
\end{equation}
The softening factor $\beta$ is chosen equal to 4, i.e., the soft material tensor is four times softer than $\mathbb{C}_{bulk}$.
We choose a dense material having $E=12.5$GPa and $\nu=0.25$.
We discretize the domain $\Omega$ using a mesh with $128\times 64$ quadrilateral elements
and we set $m=0.45$,
$\kappa_{\chi}=\kappa_{\phi}=4$,
$\gamma_{\phi}=0.02$, 
a time step increment $\Delta\tau=1.0\times 10^{-6}$, 
$\phi_0=0.5$ as initial solution, and 
a tolerance equal to 0.01.
\subsubsection{Sensitivity study of the graded-material interface parameters $\gamma_{\chi}$ }
\label{sssec:CantBeamSensitivityGammaPhi}
\begin{figure}[h!]
\centering
  \subfloat[$\gamma_{\chi}=0.001$\label{GammaChi001}]
	{
     \includegraphics[width=0.247\textwidth]{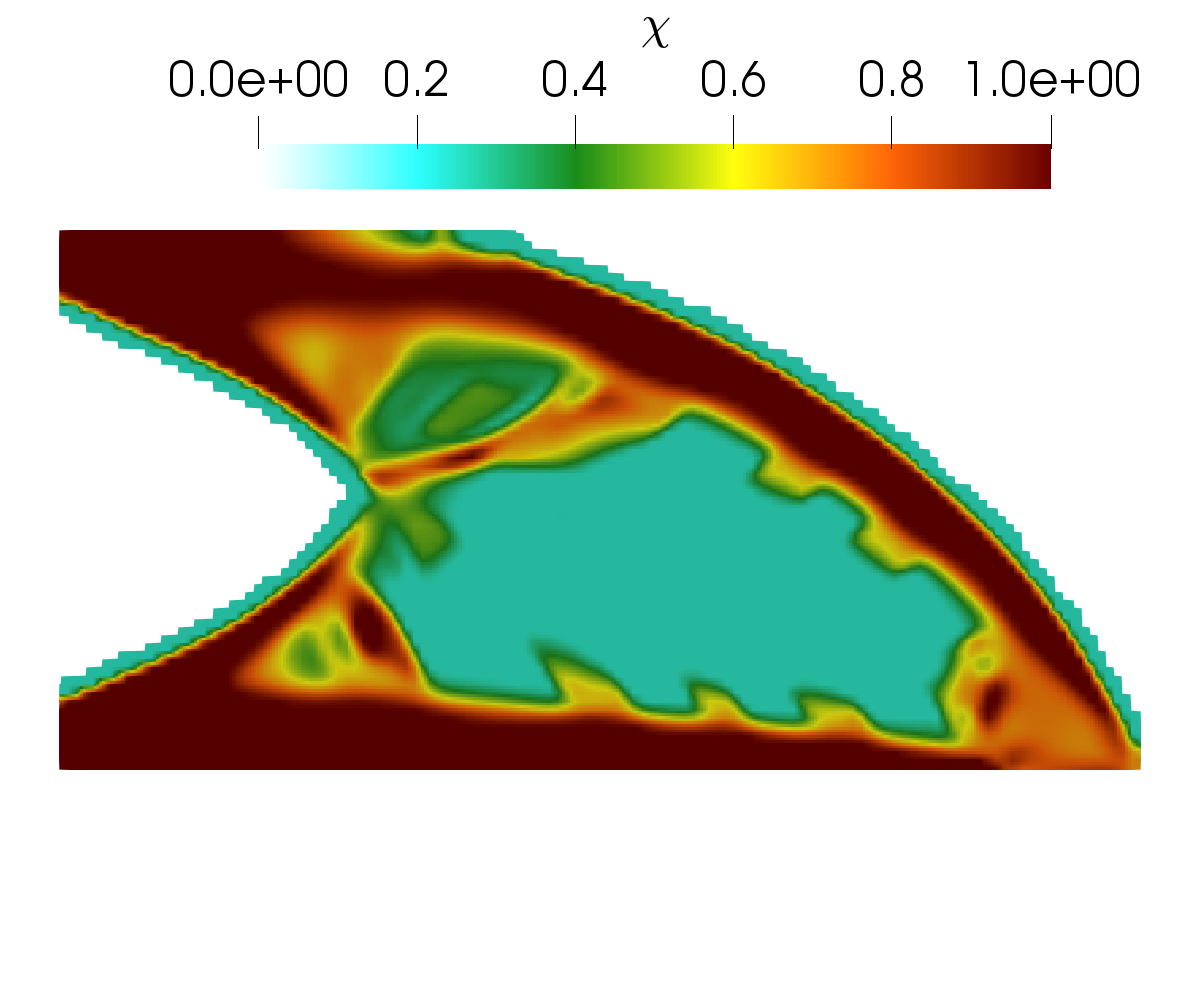}
  }
  \subfloat[$\gamma_{\chi}=0.005$\label{GammaChi005}]
	{
     \includegraphics[width=0.247\textwidth]{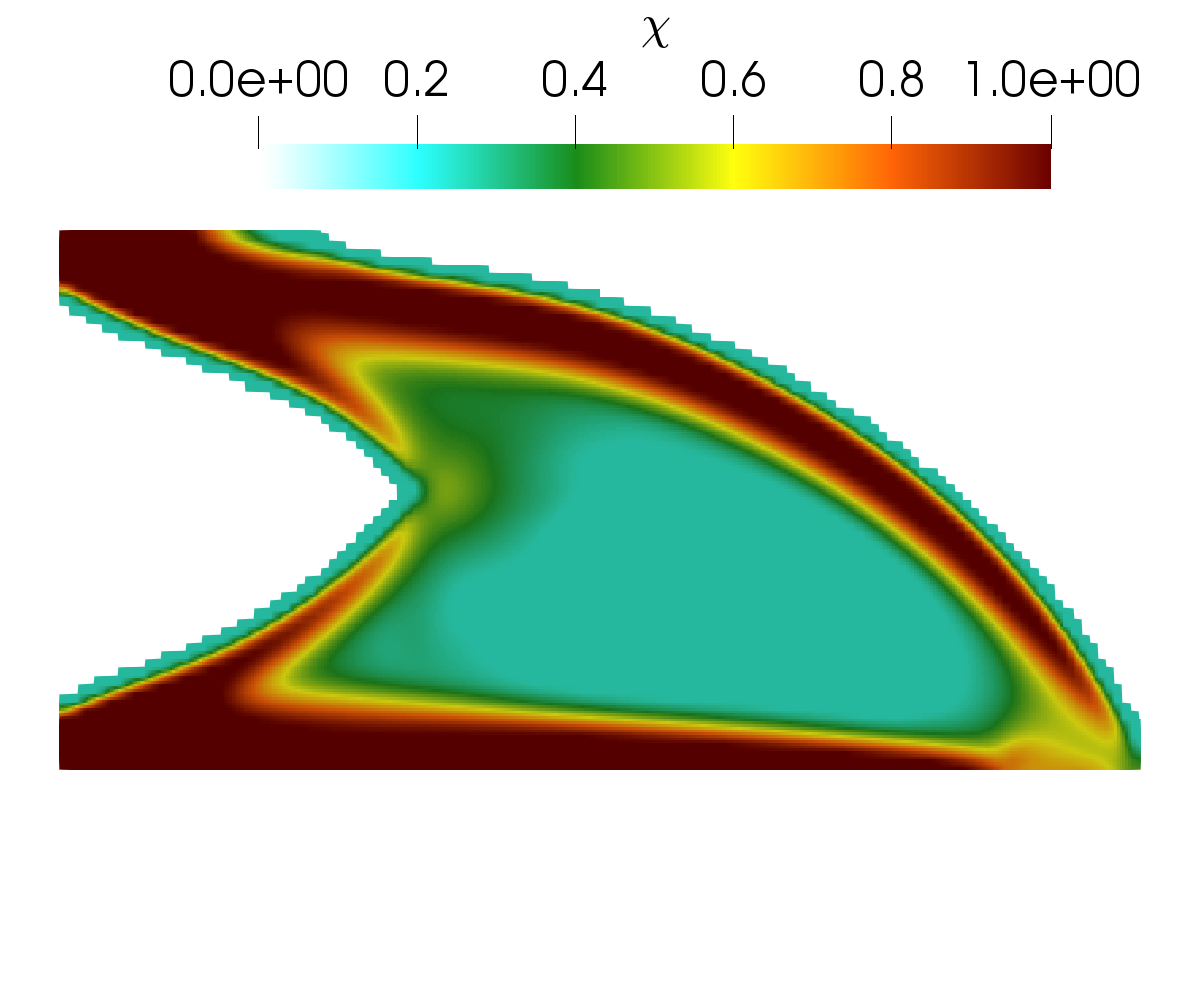}
  }
  \\
  \subfloat[$\gamma_{\chi}=0.010$\label{GammaChi010}]
	{
     \includegraphics[width=0.247\textwidth]{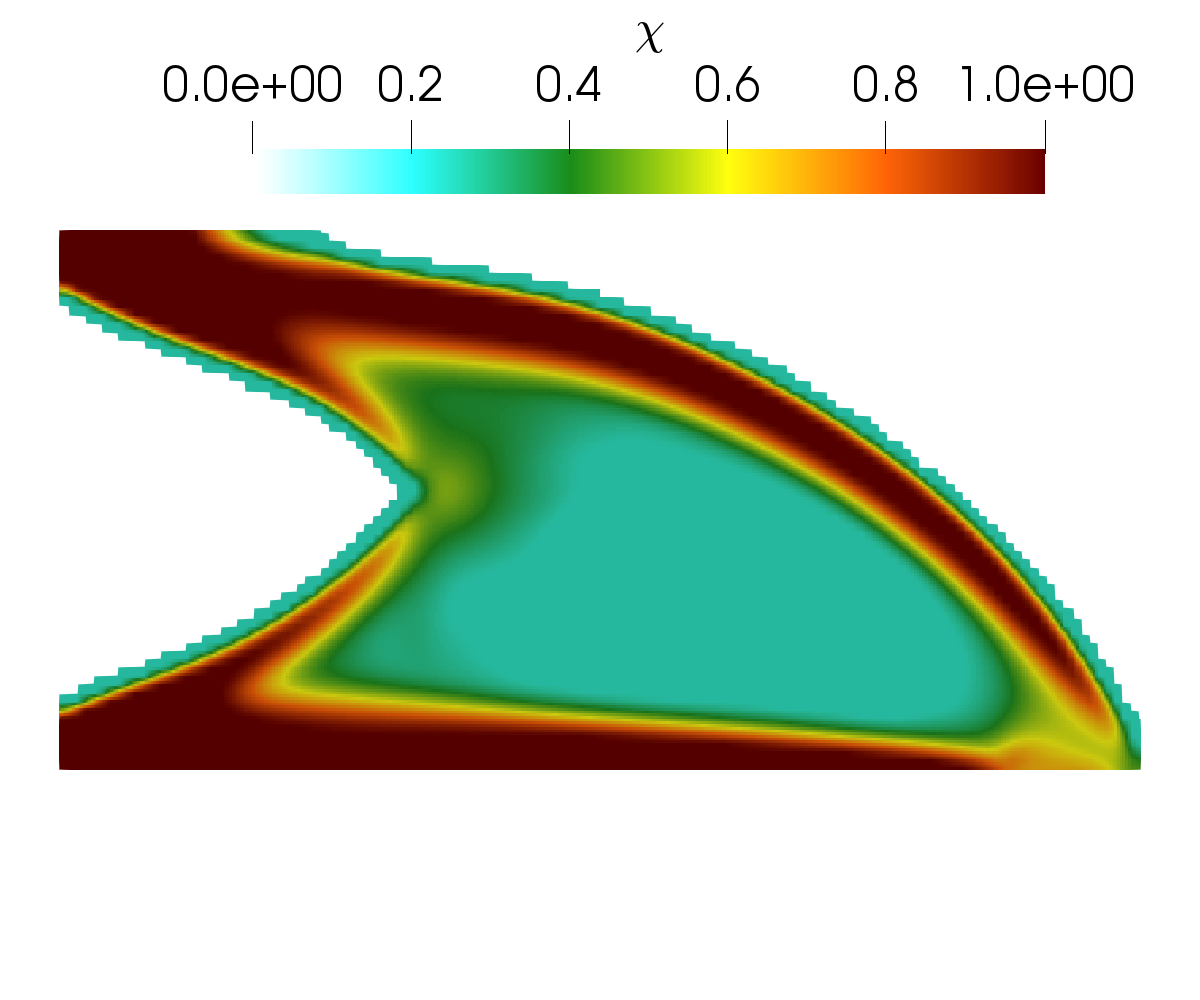}
  }
  \subfloat[$\gamma_{\chi}=0.020$\label{GammaChi020}]
	{
     \includegraphics[width=0.247\textwidth]{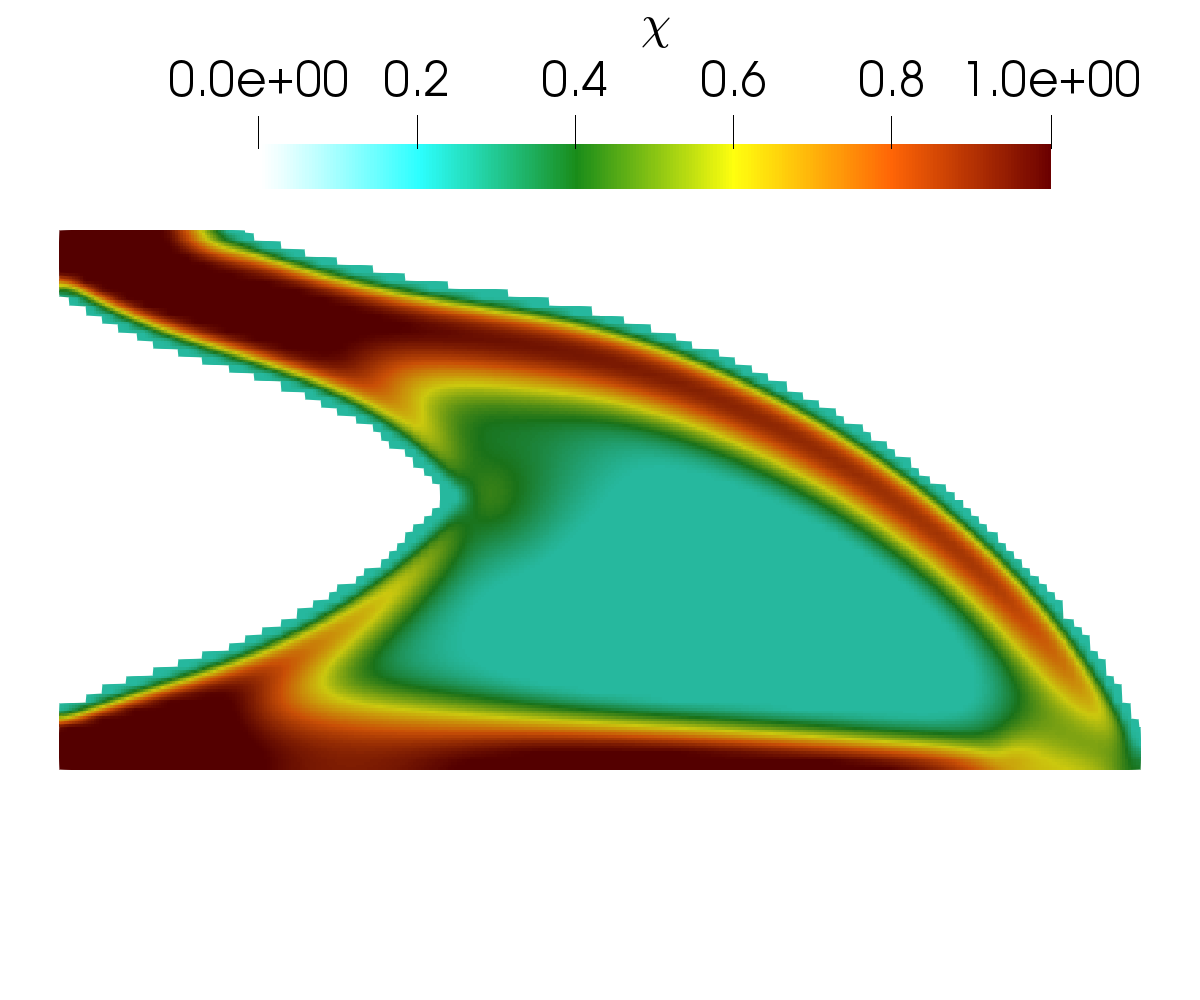}
  }
    \\
  \subfloat[$\gamma_{\chi}=0.050$\label{GammaChi050}]
	{
     \includegraphics[width=0.247\textwidth]{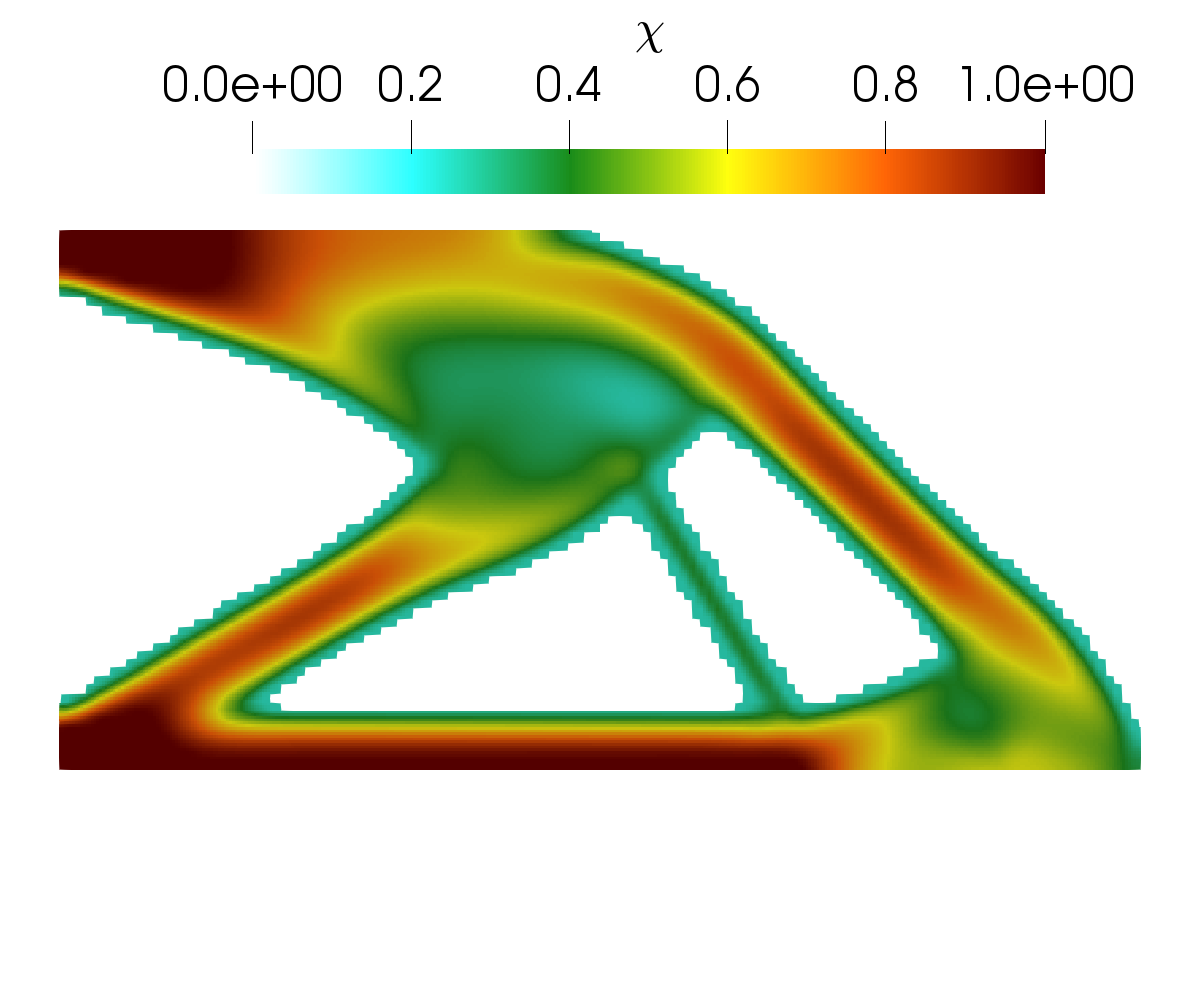}
  }
  \subfloat[$\gamma_{\chi}=0.100$\label{GammaChi1}]
	{
     \includegraphics[width=0.247\textwidth]{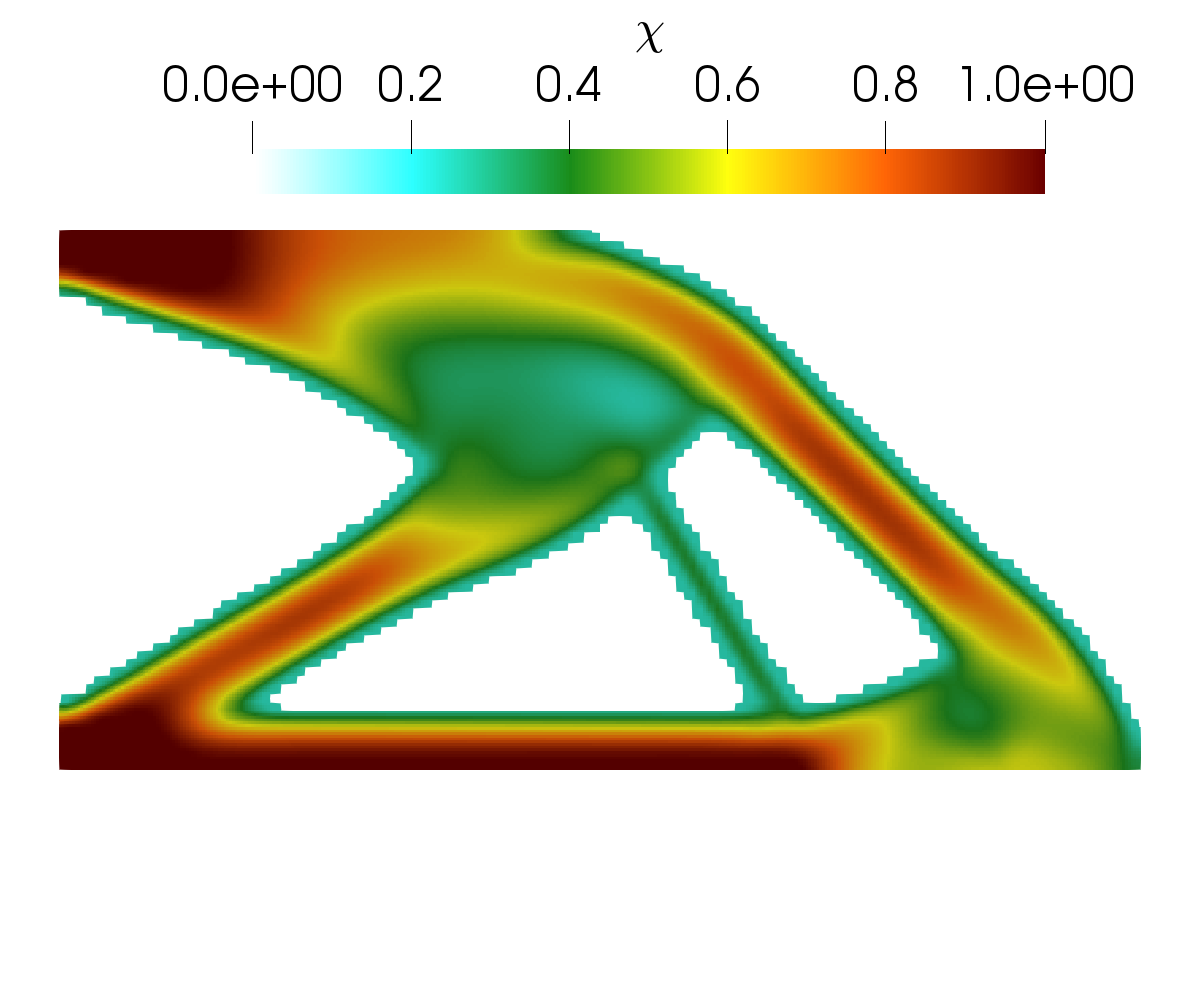}
  }
    \caption{Cantilever beam: Sensitivity study of the multi-material interface parameters $\gamma_{\chi}$. If $\gamma_{\chi}\leq \gamma_{\phi}=0.02$ different values of the graded-material interface parameters do not affect too much the final solution, which presents a wide region of soft material filling the voids which, instead, characterize the single-material solution (see~\autoref{DenseMaterial}); whereas, for $\gamma_{\chi}\geq \gamma_{\phi}$, the final solution of the graded-material case tends to a single-material configuration with multiple holes in the structure. Finally, if we choose $\gamma_{\chi}$ too small with respect to our element size, the solution does not converge anymore as \mbox{in~\autoref{GammaChi001}.}\label{sensStudy} }
\end{figure}

\begin{figure}[h!]
\centering

     \includegraphics[width=0.247\textwidth]{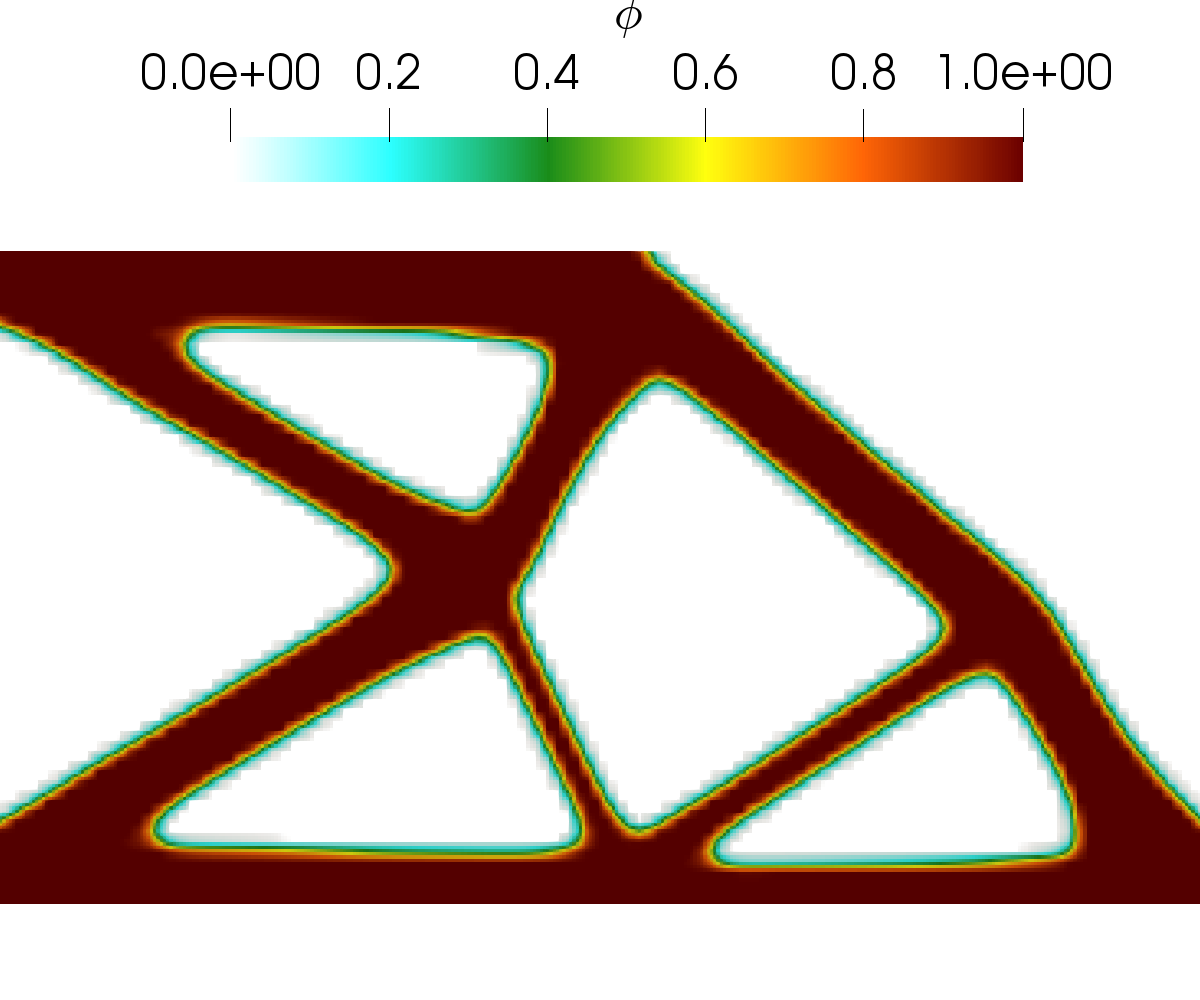}

    \caption{Cantilever beam: Sensitivity study of the graded-material interface parameters $\gamma_{\chi}$. Single-material optimized structure. \label{DenseMaterial}}
\end{figure}

In this first sensitivity study, we investigate the different topologies obtained by varying $\gamma_{\chi}$ between $0.001$ and $0.1$, as reported in~\autoref{sensStudy}.
The results show that the optimal multi-material distribution is very different from the single-material optimized topology depicted in~\autoref{DenseMaterial}, for values of $\gamma_{\chi}$ smaller than $\gamma_{\phi}$.
In fact, in this case, the voids present in the single-material structure are replaced by areas of soft material.
Contrary, if $\gamma_{\chi}$ is chosen to be bigger than $\gamma_{\phi}$ the solution presents void regions similarly to the single-material case.
Finally, we observe that, as expected, when the thickness of the diffuse interface is too small compared to the element size, the solution does not converge anymore~\mbox{(\autoref{GammaChi001})}.
\autoref{table:SensitivityGammaChi} reports the values of the compliance and of the material fraction index $m_{\chi}$ for different values of $\gamma_{\chi}$.
From the values in \autoref{table:SensitivityGammaChi}, we can see that employing a softer material will increase the compliance of the body, leading at the same time to lighter structures compared to the homogeneous material case.
\begin{table}
\centering
\caption{Cantilever beam: Sensitivity study of the graded-material interface parameters $\gamma_{\chi}$. Compliance and material index values for different choices of $\gamma_{\chi}$.}
\label{table:SensitivityGammaChi}
\begin{tabular}{llll}
\hline\noalign{\smallskip}
$\gamma_{\chi}$ & compliance& $m_{\chi}$ & convergence \\
\noalign{\smallskip}\hline\noalign{\smallskip}
$0.001$ & $105.3$ & $0.380$ & NO\\
$0.005$ & $122.9$ & $0.265$ & YES \\
$0.01$ & $133.0$ & $0.245$ & YES \\
$0.02$ & $141.9$ & $0.230$ & YES \\
$0.05$ & $154.0$ & $0.225$ & YES \\
$0.1$ & $165.4$ & $0.201$ & YES \\
\noalign{\smallskip}\hline
full dense material& $52.3$ & $m_{\phi}=0.45$ & YES \\
\noalign{\smallskip}\hline
\end{tabular}
\end{table}
\subsubsection{Sensitivity study of the slenderness of the structure}
\label{sssec:CantBeamSlanderness}
\begin{figure}[h!]
\centering
\subfloat[$s=1$]
	{
     \includegraphics[width=0.247\textwidth]{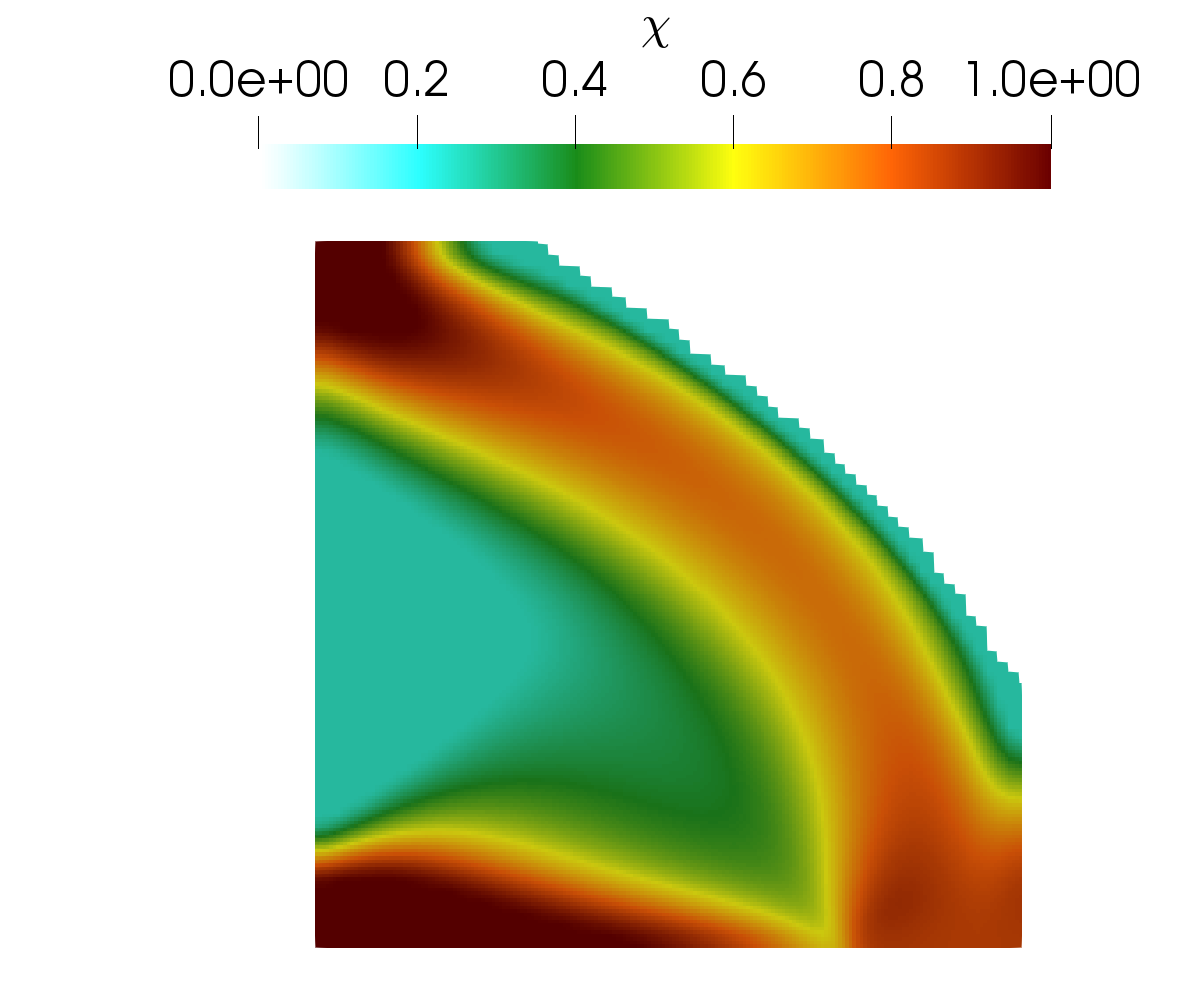}
  }
\subfloat[$s=2$]
	{
     \includegraphics[width=0.247\textwidth]{CantileverBeam_BW_Chi_M09_GammaChiNew_0020_mod.png}
  }
  \\
 \subfloat[$s=4$]
	{
     \includegraphics[width=0.247\textwidth]{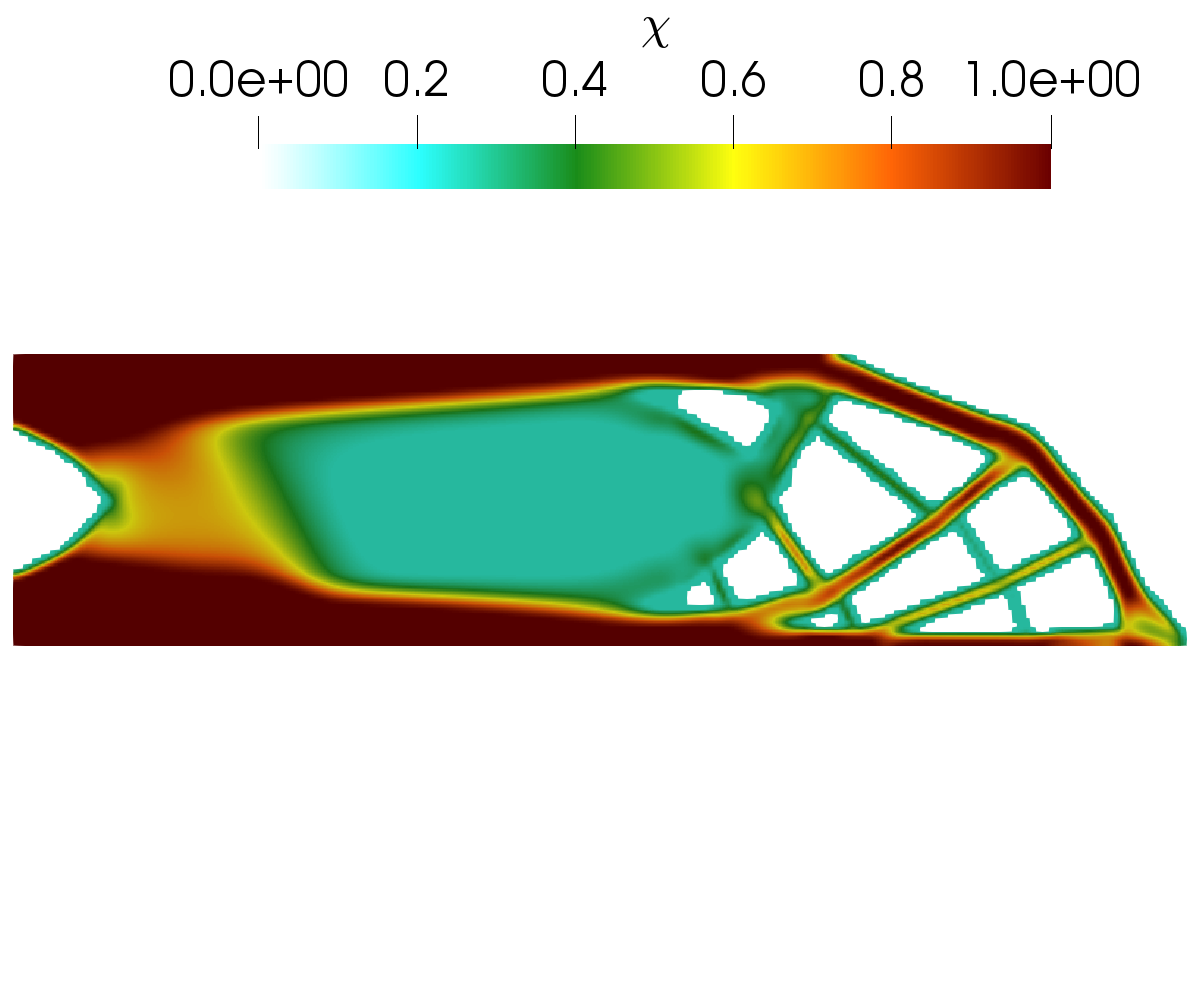}
  }
      \caption{Cantilever beam: Sensitivity study of the slanderness of the structure $s$. Varying the slanderness ratio $s=a/b$ we obtain optimized structures characterized by an outer frame of stiff material filled with regions of soft materials.\label{sensStudySlanderness} }
\end{figure}

On the cantilever beam, we perform a second sensitivity study varying the slenderness ratio $s=a/b$ (i.e., the ratio between the length and the height of the cantilever beam), for a fixed value of the graded-material interface parameter ($\gamma_{\chi}=0.02$).
\autoref{sensStudySlanderness} shows the final topologies for three different slenderness ratios, where all the resulting structures are characterized by internal regions of softer material and external support of stiffer material.
\subsection{Simply-supported beam}
\label{ssec:bridge}
\begin{figure}
\centering
\includegraphics[width=0.497\textwidth]{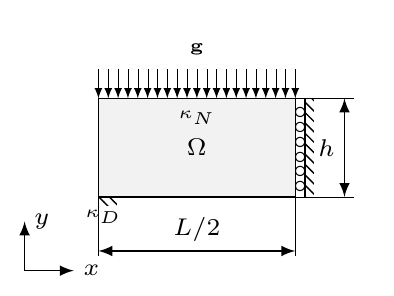}
\caption{Simply-supported beam: Initial configuration and problem domain.\label{fig:BridgeProblem}}
\end{figure}

In this second example, we choose Acrylonitrile Butadiene Styrene (ABS), which is a common thermoplastic polymer widely used in 3D printing applications, as material to obtain an optimized simply-supported beam structure.
The problem is symmetric and thus we decide to solve only half of the domain as depicted in~\autoref{fig:BridgeProblem}, where $h=50$mm an $L/2=100$mm,
with a distributed external load $\mathbf{g}$ equal to $(0,-50)$N/mm applied on the top edge of the structure.
The Young modulus and the Poisson coefficient of ABS plastic are $2.3$GPa and $0.35$, respectively.
We set $m=0.4$,
$\kappa_{\phi}=\kappa_{\chi}=1$
and $\gamma_{\phi}=\gamma_{\chi}=0.01$,
while we choose a pseudo-time step $\Delta\tau=1.0\times 10^{-6}$ and an initial solution $\phi_0=0.5$.
\subsubsection{Sensitivity study of the softening factor $\beta$}
\label{sssec:bridgeSoftening}
\begin{figure}[ht]

\centering
  \subfloat[$\beta=1$\label{Beta1}]
	{
     \includegraphics[width=0.247\textwidth]{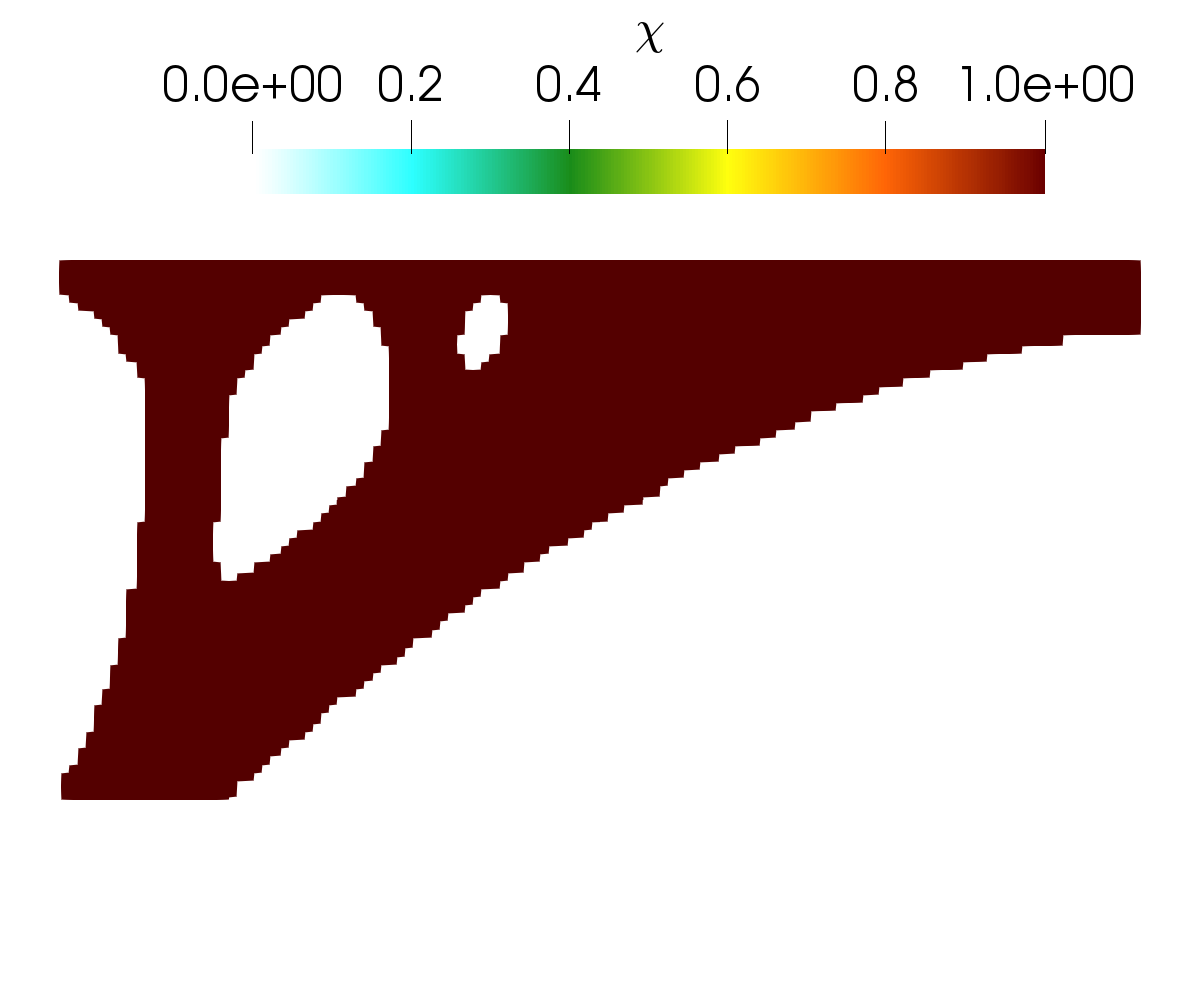}
  }
\subfloat[$\beta=2$\label{Beta2}]
	{
     \includegraphics[width=0.247\textwidth]{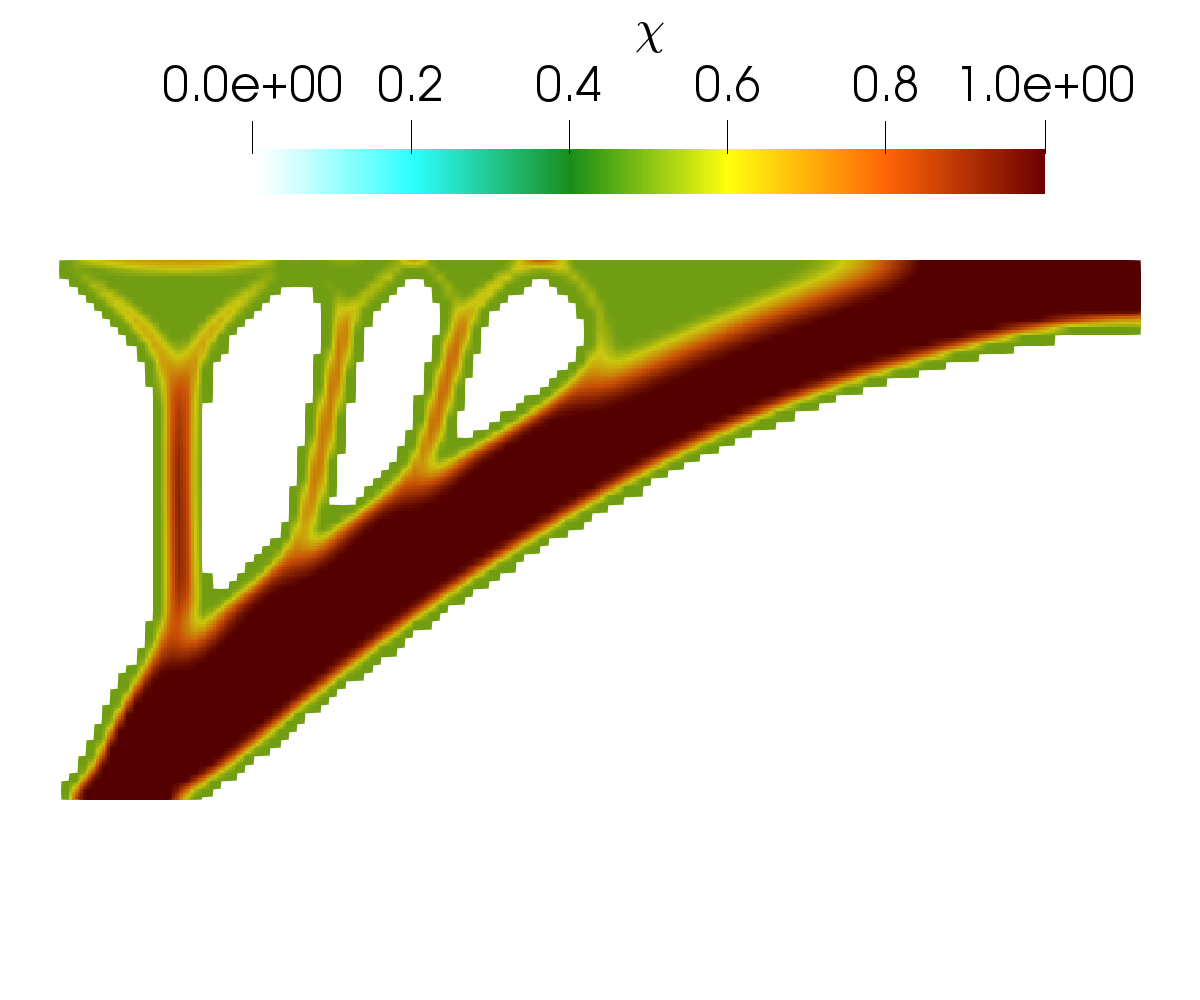}
  }
  \\
  \subfloat[$\beta=3$\label{Beta3}]
	{
     \includegraphics[width=0.247\textwidth]{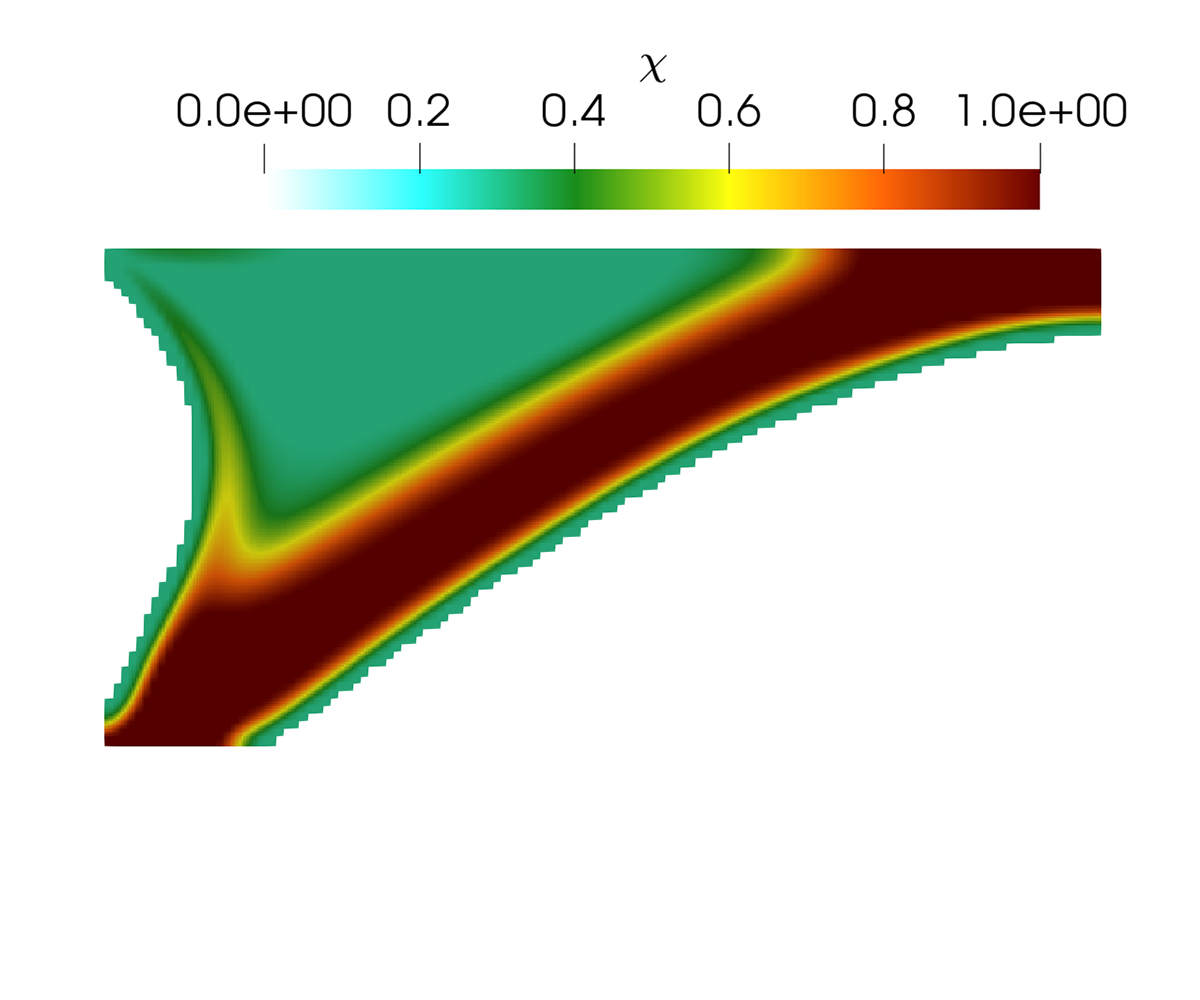}
  }
    \subfloat[$\beta=4$\label{Beta4}]
	{
     \includegraphics[width=0.247\textwidth]{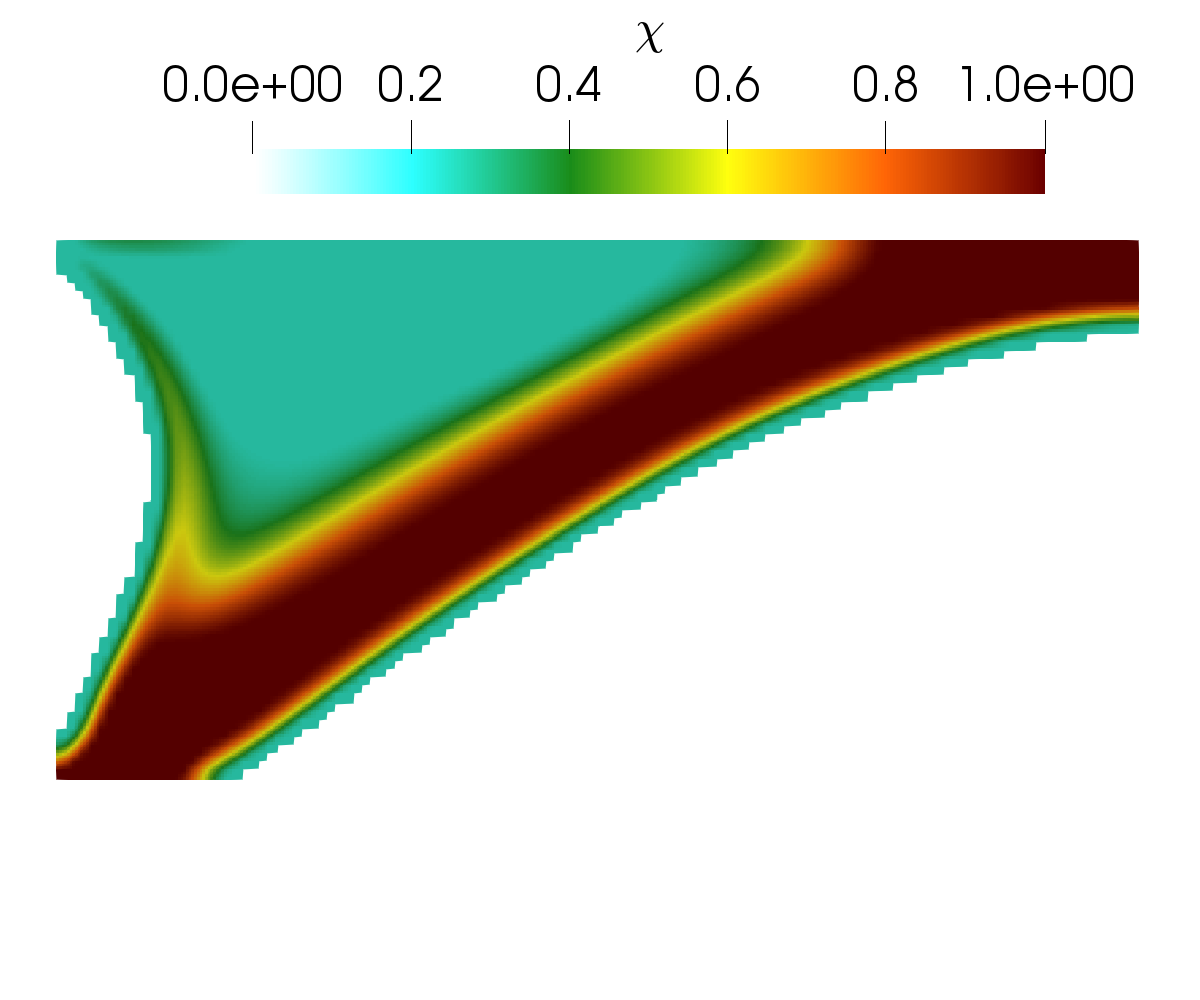}
  }

\caption{Simply-supported beam: Sensitivity study of the softening factor $\beta$. Increasing the values of the softening factor, i.e., employing a softer material, the optimized structure does not present anymore the typical holes resulting from a single-material optimization~\ref{Beta1}. Voids are now replaced by a region of soft material. \label{fig:BridgeTopEv}}
\end{figure}

\autoref{fig:BridgeTopEv} presents the results of a sensitivity analysis performed varying the softening factor $\beta$ from 1 to 4.
The resulting optimized structures show that, introducing \textit{grey-scale} regions in the structure, the optimal design is modified, replacing the typical voids of SIMP approach with areas of soft material.
Again we observe in~\autoref{table:SensitivityBeta} that introducing a soft material within the algorithm leads to structure with a smaller material index but higher compliance. 
The values of the softening factor strongly influences the final results and give us the possibility to obtain intermediate structure such as the one in~\autoref{Beta2}. Moreover, we can notice that for high values of $\beta$ the results are very similar to each other (see~\autoref{Beta3}and~\autoref{Beta4}), thus the higher values of this parameter would depend only on the technological boundaries of the AM process.
\begin{table}
\centering
\caption{Simply-supported beam: Sensitivity study of the softening factor $\beta$. Compliance and material index values for different choices of $\beta$.}
\label{table:SensitivityBeta}
\begin{tabular}{llll}
\noalign{\smallskip}\hline
$\beta$ & compliance & $m_{\chi}$ & convergence \\ 
\noalign{\smallskip}\hline\noalign{\smallskip}
1 & $20.5$ & $0.40$ & YES\\
2 & $37.3$ & $0.32$ & YES \\
3 & $46.4$ & $0.24$ & YES \\
4 & $58.6$ & $0.18$ & YES \\
\noalign{\smallskip}\hline
\end{tabular}
\end{table}
\subsubsection{Sensitivity study of the distributed load}
\label{sssec:BridgeLoad}
\begin{figure}[ht]
\centering
  \subfloat[$\mathbf{g}\times2$]
	{
     \includegraphics[width=0.247\textwidth]{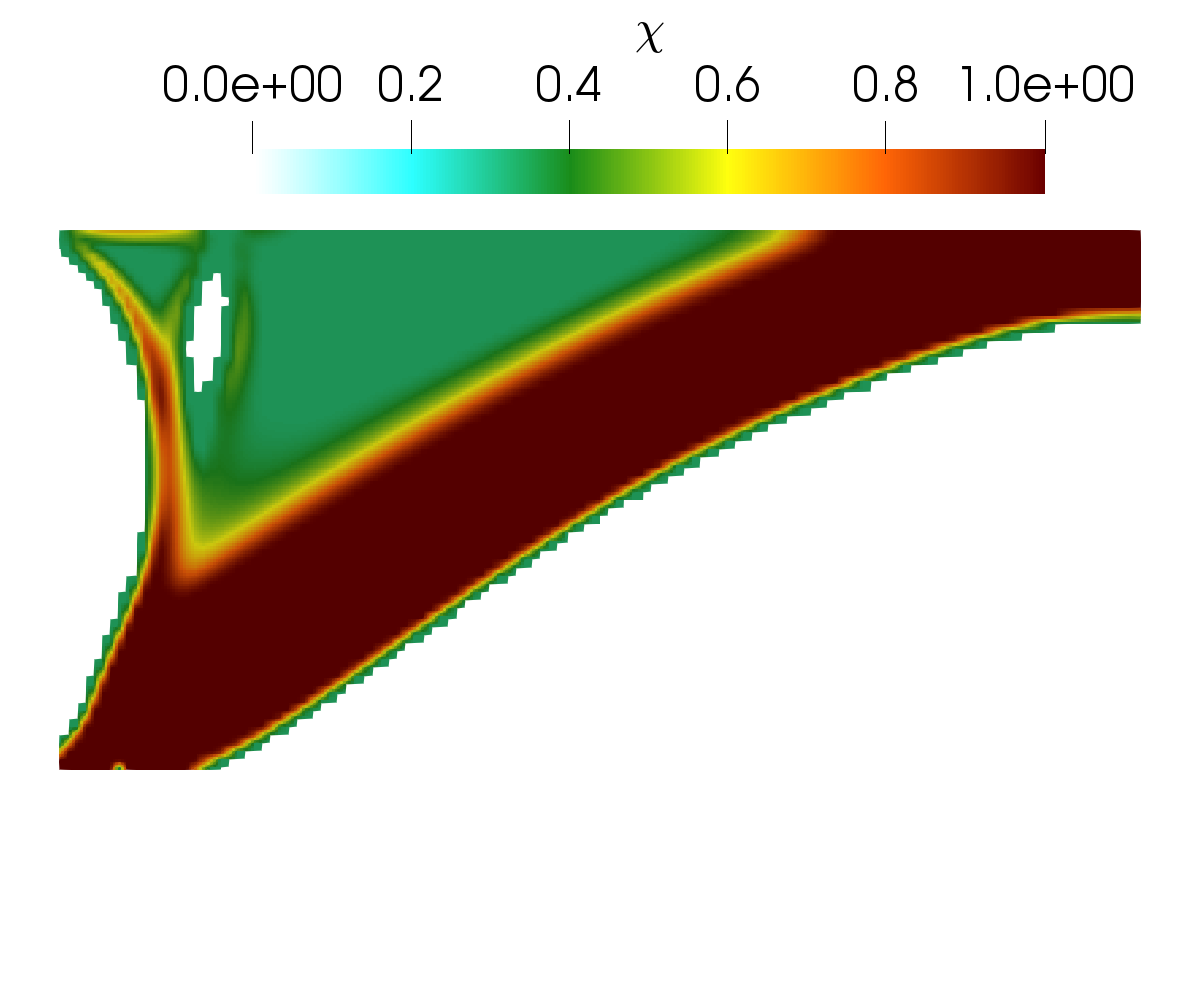}
  }
\subfloat[$\mathbf{g}\times3$ \label{fig:BridgeTopg3}]
	{
     \includegraphics[width=0.247\textwidth]{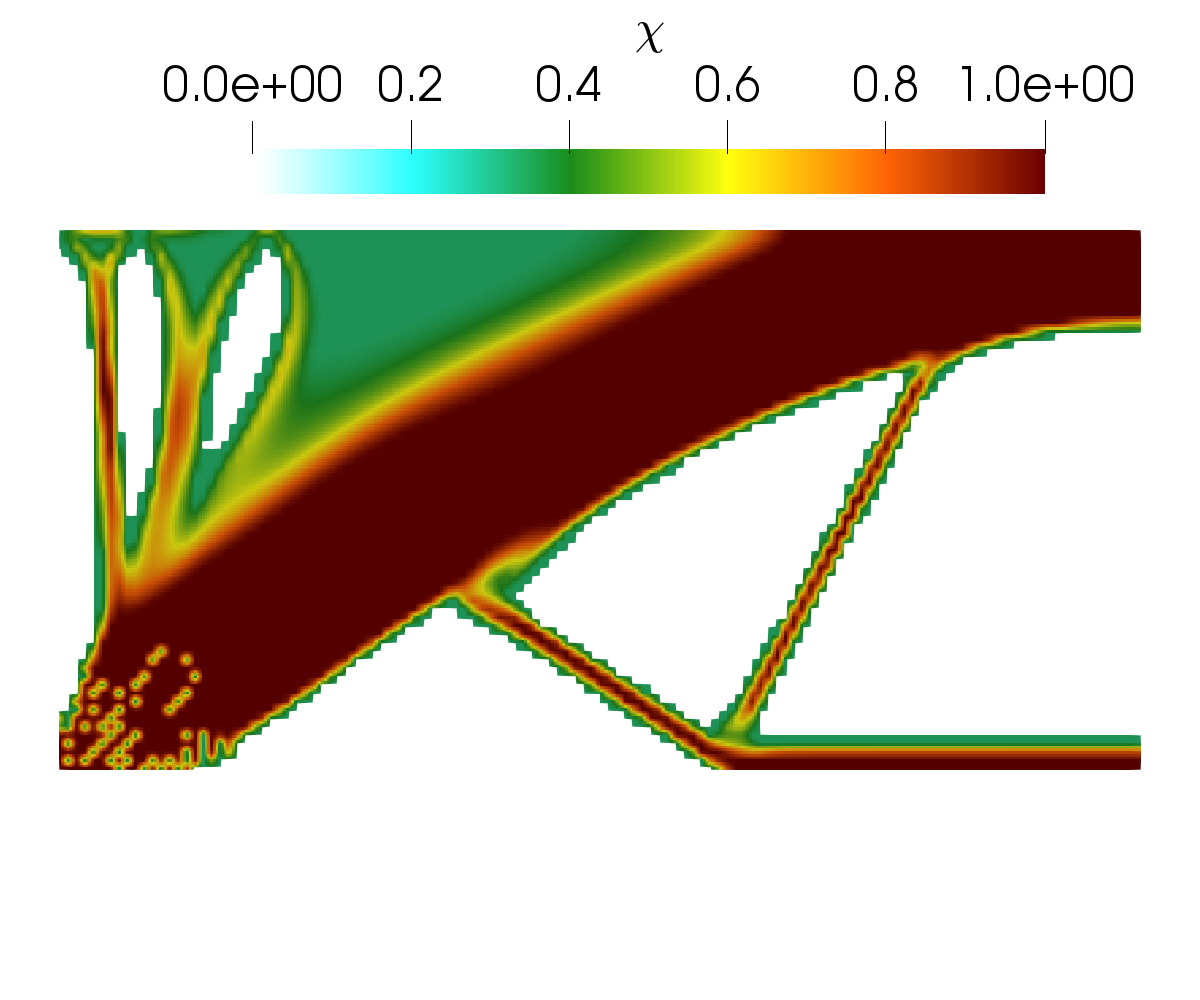}
  }

\caption{Simply-supported beam: Sensitivity study of the distributed load $\mathbf{g}$. Increasing the load on the upper edge of the structure we observe an increment in the region occupied by full dense material and the presence of holes and columnar structures. \label{fig:BridgeTopLoadInc}}
\end{figure}

On the simply-supported beam, we conduct a second sensitivity study fixing $\beta=3$ and increasing the distributed load $\mathbf{g}$ by a factor of 2 and 3, respectively. 
The resulting structures are reported in~\autoref{fig:BridgeTopLoadInc}.
As we expected employing a heavier load reduces the areas of soft material, increasing at the same time the  number of columnar structures in the final topology.
We want to remark here that the structure of~\autoref{fig:BridgeTopg3} did not converge even after 1000 iterations. 
Since the mesh is not modified, this behavior is due to the choice of the stiffer material, which is in this case too soft for such a heavy load.
\subsection{From numerics to AM products}\label{ssec:NumToAM}
\begin{figure}
\centering
\includegraphics[width=0.457\textwidth]{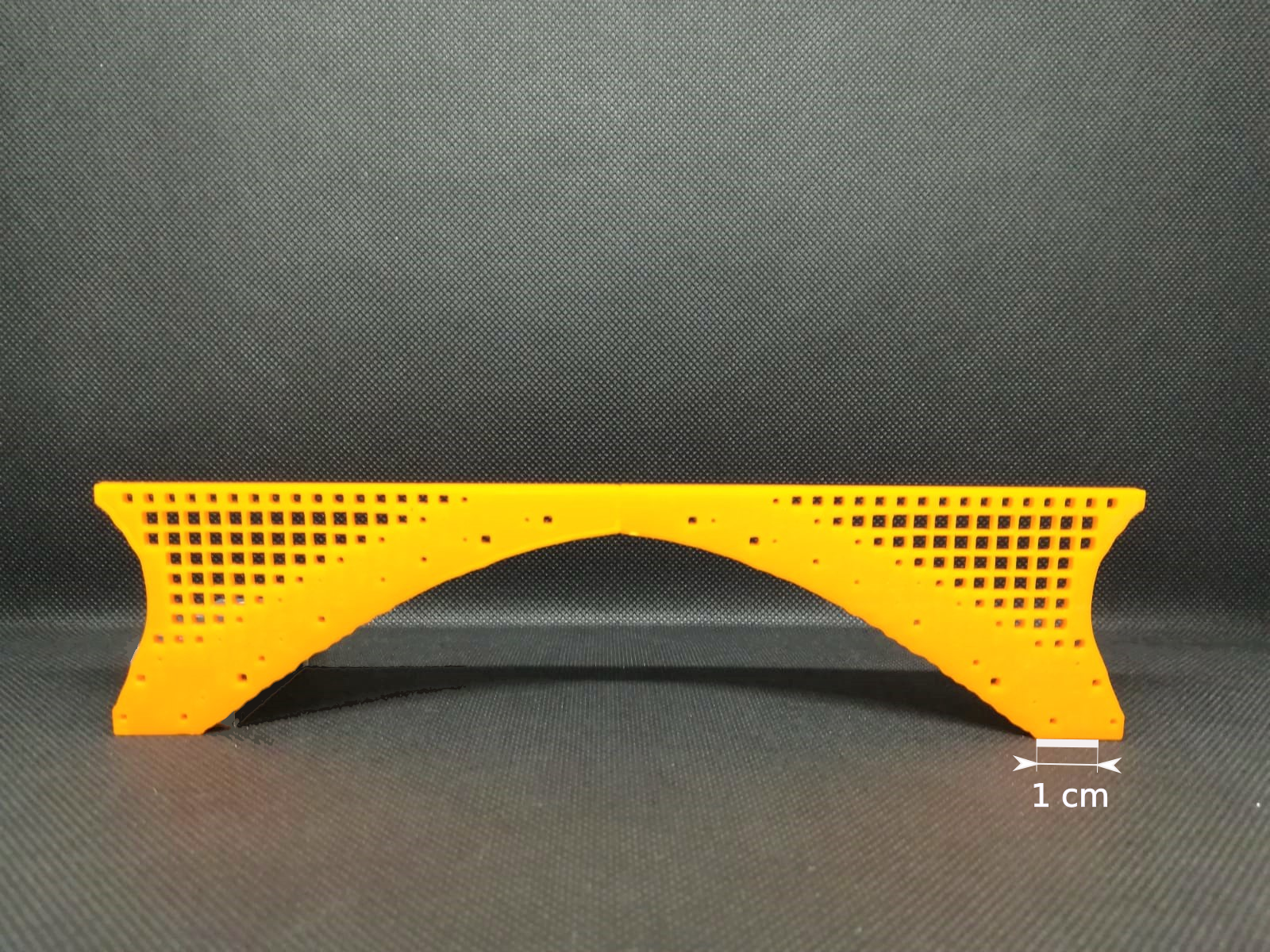}
\caption{Simply-supported beam with a distributed load of \autoref{Beta4} printed using FDM 3D printer.}\label{Bridge3Dprinted}
\end{figure}

In order to demonstrate that the numerical results presented so far can be actually realized in practice, we decided to print a sample of the optimized structure depicted in~\autoref{Beta4}. To obtain this result we used the Fused Deposition Modeling (FDM) 3D printer present in our PROTOtyping LABoratory (PROTOLAB) at the University of Pavia. This sample is made in ABS plastic material and realized extruding by an offset of 5mm the structure of ~\autoref{Beta4}. The values of $\chi$ are then mapped onto a manufacturing grid where the dimension of each cell depends on technological constraints given by the machine. Finally, by means of a simple boolean operation, we subtract a quadratic region proportional to the average value of $\chi$ in each cell. 
Once the corresponding CAD model is completed we can directly print it, obtaining the structure of~\autoref{Bridge3Dprinted}. Since, the assumption of linear elasticity is not valid for a soft material such as plastic (metal alloys would be a more appropriate choice in this case) we leave the experimental validation of the proposed numerical algorithm to future research.
\section{Conclusions} \label{sec:conclusions}
In the present work, we have introduced a novel phase-field topology optimization algorithm based on a graded material definition.
\par
The numerical results show that the additional control parameter $\chi$, introduced in our phase-field formulation, allows increasing the number of possible optimal designs delivered by the topology optimization process.
\par
In particular, we have introduced the possibility to control the distribution of the material density within our structure in a continuous fashion.
Such a feature can be in many cases highly desirable, in particular if we consider additive manufacturing applications.
\par
Moreover the algorithm allows to easily control the number of regions with graded material distribution, delivering results which can be in between a fully \textit{black-and-white} approach and a purely graded-material distribution.
\par
In the near future we aim at investigating mechanical properties of 3D printed structures designed using the graded-material phase-field topology optimization algorithm. 

\begin{acknowledgements}
This work was partially supported by Regione Lombardia through the project "TPro.SL - Tech Profiles for Smart Living" (No. 379384) within the Smart Living program, and through the project "MADE4LO - Metal ADditivE for LOmbardy" (No. 240963) within the POR FESR 2014-2020 program.
MC and AR have been partially supported by Fondazione Cariplo - Regione Lombardia through the project ``Verso nuovi strumenti di simulazione super veloci ed accurati basati sull'analisi isogeometrica'', within the program RST - rafforzamento. This research has been performed in the framework of the project Fondazione Cariplo-Regione Lombardia MEGAsTAR ``Matematica d'Eccellenza in biologia ed ingegneria come
acceleratore di una nuova strateGia per l'ATtRattivit\`a dell'ateneo pavese''. The present paper also benefits from the support of the GNAMPA (Gruppo Nazionale per l'Analisi Matematica, la Probabilit\`a e le loro Applicazioni) of INdAM (Istituto Nazionale di Alta Matematica) for ER.
A grateful acknowledgment goes to Dr. Ing. Gianluca Alaimo for his support and precious suggestions on additive manufacturing technology.
\end{acknowledgements}

\bibliographystyle{unsrtnat}
{{\footnotesize \bibliography{bibliography}}

\end{document}